\documentclass[10pt]{article}
\usepackage[margin=1in]{geometry}
\usepackage{comment}
\usepackage[T1]{fontenc}
\usepackage{lmodern}
\usepackage{todonotes}
\usepackage{ulem}
\usepackage{bm}

\usepackage[colorlinks=true,linkcolor=red,citecolor=blue]{hyperref}%
\usepackage{cite}
\usepackage{url}            
\usepackage{booktabs}       
\usepackage{amsfonts} 
\usepackage{nicefrac}       
\usepackage{microtype} 
\usepackage{wrapfig}
\usepackage{enumerate}
\usepackage{graphicx}
\usepackage{amsmath,amssymb}
\usepackage{bigstrut}
\usepackage{float}
\usepackage{caption}
\usepackage{algorithm}
\usepackage{mathtools}
\usepackage{subfigure}
\usepackage{booktabs}
\usepackage{authblk}
\usepackage{algpseudocode}
\usepackage{dsfont}
\usepackage{bbm}
\usepackage{comment}

\newtheorem{theorem}{Theorem}

\newtheorem{remark}{Remark}
\newtheorem{lemma}{Lemma}
\newtheorem{example}{Example}

\newtheorem{asm}{Assumption}

\def\R{\mathbb{R}}
\def\ones{\mathbf{1}}
\def\xh{\hat{\bfx}}
\def\nt{\bfe}
\def\sm{\gamma}

\newcommand{\Q}{Q}
\newcommand{\C}{C}

\newcommand{\EE}[1]{\mathbb{E}\left[ #1 \right]}

\newcommand{\s}{\star}
\newcommand{\bfx}{\mathbf{x}}
\newcommand{\bfe}{\mathbf{e}}
\newcommand{\bfy}{\mathbf{y}}
\newcommand{\bfv}{\mathbf{v}}

\newcommand{\bfz}{\mathbf{z}}
\newcommand{\bfr}{\mathbf{r}}

\newcommand{\m}{\rho} 
\newcommand{\li}{L} 

\newcommand{\bfu}{\boldsymbol{u}}

\newcommand{\clg}{\mathcal{G}}

\newcommand{\cle}{\mathcal{E}}

\def\V{V_\bfr}
\def\rmin{\bfr_{\min}}
\def\Wt{\{W(t)\}}
\def\minW{\eta}

\newcommand{\Del}{\Delta}
\renewcommand{\mod}{\mathsf{\ mod\ }}
\renewcommand{\comment}[1]{}

\newcommand{\F}{\mathcal{F}}
\newcommand{\azr}{\alpha_0}
\newcommand{\bzr}{\beta_0}

\newcommand{\md}{\middle| }
\DeclareMathOperator{\diag}{diag}

\newcommand{\sign}[1]{\operatorname{sgn}(#1)}

\newcommand{\mx}{\bar{\bfx}}
\newcommand{\mg}{\bar{g}}
\newcommand{\lp}{\left(}
\newcommand{\rp}{\right)}
\newcommand{\lb}{\left[}
\newcommand{\rb}{\right]}
\newcommand{\la}{\left\langle}
\newcommand{\ra}{\right\rangle}
\newcommand{\lnr}{\left\|}
\newcommand{\rnr}{\right\|}
\newcommand{\lc}{\left\{}
\newcommand{\rc}{\right\}}
\newcommand{\dl}{\delta}
\newcommand{\et}{\beta(t)}
\newcommand{\at}{\alpha(t)}

\newcommand{\BO}{\mathcal O}

\newcommand{\ER}[1]{E(#1)}

\newcommand{\cvxTwo}{\epsilon_1}
\newcommand{\cvxThree}{\epsilon_2}

\newcommand{\cvxFive}{\epsilon_4}
\newcommand{\cvxSix}{\epsilon_5}

\newcommand{\cOne}{\xi_1}
\newcommand{\cTwo}{\xi_2}

\newcommand{\Tzero}{T_1}
\newcommand{\Tfirst}{T_2}
\newcommand{\TStrt}{T_3}

\newcommand{\TStrtTwo}{\tau}
\newcommand{\kap}{\kappa}
\newcommand{\tLM}{\tau}
\newcommand{\Ts}{T_0}

\newcommand{\Ta}{T_4}
\newcommand{\Tb}{T_5}
\newcommand{\Tc}{T_6}
\newcommand{\Td}{T_7}
\newcommand{\Te}{T_8}

\newcommand{\Nr}[1]{\lnr #1 \rnr_{\bfr}}

\title{Distributed Optimization over Time-varying Graphs \\ with Imperfect Sharing of Information}

\author{Hadi Reisizadeh, Behrouz Touri, and Soheil Mohajer\footnote{H.\ Reisizadeh (email: hadir@umn.edu) and S.\ Mohajer (email: soheil@umn.edu) are with the University of Minnesota, and B.\ Touri (email: btouri@ucsd.edu) is with the University of California San Diego.}}
\date{} 
\begin{document}
\maketitle

\begin{abstract}
We study strongly convex distributed optimization problems where a set of agents  are interested in solving a separable optimization problem collaboratively. In this paper, we propose and study a two time-scale decentralized gradient descent algorithm for a broad class of lossy sharing of information over time-varying graphs. One time-scale fades out the (lossy) incoming information from neighboring agents, and one time-scale regulates the local loss functions' gradients. For strongly convex loss functions, with a proper choice of step-sizes, we show that the agents' estimates converge to the global optimal state at a rate of $\BO \lp T^{-1/2}\rp$. Another important contribution of this work is to provide novel tools to deal with diminishing average weights over time-varying graphs. 
\end{abstract}

\section{Introduction}
Emergence of big data analytics, modern computer architectures, storage, and data collection have led to a growing interest in the study of multi-agent networked systems. These systems arises in various applications such as sensor networks~\cite{rabbat2004distributed,kar2008distributed}, network routing~\cite{neglia2009distributed}, large scale machine learning~\cite{tsianos2012consensus}, power control~\cite{ram2009distributed}, and distributed network resource allocations~\cite{xiao2006optimal,ribeiro2010ergodic}, for which decentralized solutions offer promising results. In general and in the absence of a central entity, we are often dealing with a time-varying network of agents,   each can perform local and on-device computation. The information can be shared throughout the network via local communication between neighboring agents. This communication among agents, specially when the dimension of the data is large, accounts for a significant delay in the overall running time of the algorithm. In this paper, we study such a distributed optimization framework with lossy and imperfect information sharing, and propose and analyze an gradient-based distributed algorithm which guarantees convergence to the optimum solution, in spite of a limitation on the communication load.

\noindent\textbf{Related Works.} 
Various methods have been proposed and studied to solve distributed optimization problems in convex settings~\cite{schizas2007consensus,nedic2008distributed,nedic2010constrained,boyd2011distributed,jakovetic2014fast,yuan2016convergence,boyd2006randomized,dimakis2010gossip,duchi2011dual}, strongly convex settings \cite{yuan2016convergence,qu2017harnessing,xi2017dextra}, and non-convex settings~\cite{tatarenko2017non,zeng2018nonconvex}. For the convex objective functions, a sub-gradient method with a fixed step-size is proposed over time-varying graphs in~\cite{nedic2009distributed}. It is shown that the objective cost function reduces at rates of $\BO\lp T^{-1}\rp$ until it reaches a neighbor of a minimizer of the original problem. To achieve exact convergence to a minimizer, various diminishing step-size sub-gradient methods have been proposed and studied~\cite{nedic2010constrained,nedic2013distributed,nedic2014distributed,jakovetic2014fast,tatarenko2017non,zeng2018nonconvex}. Considering convex loss functions that are Lipschitz continuous and have bounded gradients, a subgradient-push algorithm is proposed in~\cite{nedic2014distributed}. There it is shown that the objective cost function convergences at the rate of $\BO \lp T^{-1/2}\ln T\rp$ over uniformly strongly connected, directed time-varying graphs. Under the same assumption and strong-convexity for loss functions, a better rate $\BO\lp T^{-1}\ln T \rp$ for the objective loss function plus squared consensus residual is shown in \cite{nedic2013distributed}.

Almost all the aforementioned works on this domain, consider distributed optimization with perfect sharing of information, i.e., the agents are allowed to communicate real-valued vectors perfectly over perfect communication channels. However, exchanging exact information among nodes initiates a massive communication overhead on the system that considerably slows down the convergence rate of these algorithms in real world applications. Thus, it is reasonable to assume that each agent has access to a lossy version of neighboring agents' information. 

To address lossy/noisy sharing of information, a (fixed steps-size) decentralized gradient descent method is proposed in~\cite{reisizadeh2019exact}. Assuming fixed communication network and strongly convex local cost functions, there it is shown that for a given iteration $T$, the algorithm's parameters (depending on $T$) can be chosen such that the local estimate of each agent at iteration $T$ is (roughly) within $c\lp T^{-1/2 + \epsilon}\rp$-distance of the global optimal solution for some $c>0$ and any $\epsilon>0$. Furthermore, the result holds for a termination time $T$ which is required to satisfy $T\geq T_{\min}$, where $T_{\min}$ depends on $\epsilon$ as well as non-local parameters of the underlying fixed graph. Specifically, as $\epsilon$ goes to zero, $T_{\min}$ diverges to infinity. In a closely related recent work~\cite{vasconcelos2021improved}, a two time-scale gradient descent algorithm has been presented for strongly convex loss functions. Assuming a \textit{fixed} topology for the underlying network, uniform weighting of the local cost functions, and a specific scheme for lossy sharing of information, it is shown that the expected objective loss function achieves a rate of $\BO(T^{-1/2}(\ln T)^2)$. In another related work~\cite{srivastava2011distributed}, a two-time-scale gradient descent algorithm was presented for distributed constrained and convex optimization problems over an i.i.d.\ communication graph with noisy communication links, and sub-gradient errors. It is shown that under certain conditions on the i.i.d.\ communication graph and proper choices of time-scale parameters the proposed dynamics results in almost sure convergence of local states to the optimal point.

\noindent\textbf{Contributions.} In this work, we study distributed optimization problem for a \textit{broad} class of lossy/noisy sharing of information over \textit{time-varying} communication networks. The learning method relies only on local computations and received imperfect information from neighbor agents. We show that a two-time scale gradient descent algorithm with a proper choice of parameters, reaches the global optima (in $L_2$ and hence, in probability) for \textit{every} agent with a rate of $\BO \lp T^{-1/2}\rp$.

In addition, in the existing works on distributed optimization \cite{nedic2008distributed,nedic2009distributed,nedic2013distributed,jakovetic2014fast,yuan2016convergence, tatarenko2017non,xi2017dextra,reisizadeh2019exact,vasconcelos2021improved} (with perfect or imperfect sharing of information), either the underlying communication network is assumed to be fixed, or the non-zero elements of averaging weights are assumed to be uniformly bounded away from zero. In our proposed method, however, the weights are not uniformly bounded away from zero and they are evolving over an underlying time-varying communication networks. One of the key contributions of this paper is to develop tools and techniques to deal with diminishing averaging weights for distributed optimization over time-varying networks.

\noindent\textbf{Outline of the Paper.} In Section~\ref{sec:prob-def-main-res}, we introduce the distributed optimization problem of interest, state the main result of this work, and discuss some of its immediate implications. To support our theoretical analysis, we provide some simulation results in Section~\ref{sec:Exp-result}. The rest of the paper is devoted to the proof of the main result. In order to do so, first we provide some useful tools and results in form of intermediate lemmas in Section~\ref{sec:aux-lemmas}, whose proofs are presented in Appendix. In Section~\ref{sec:proof}, we present the proof of the main result by combining (i) the study of the deviation of each agent's state from the average state over time-varying networks, and (ii) the study of the deviation of the average state from the optimal state. Finally, we conclude the paper in Section~\ref{sec:con}.

\noindent\textbf{Notation.} 
Throughout this paper, we denote the set of integers $\{1,2,\dots,n\}$ by $[n]$ and the set of non-negative real numbers by $\mathbb{R}^+$. In this paper, we are dealing with $n$ agents that are minimizing a function in $\R^d$. For notational convenience, throughout this paper, we assume that the underlying functions are acting on \textbf{row} vectors, and hence, we view vectors in  $\R^{1\times d}=\R^d$ as row vectors. The rest of the vectors, i.e., the vectors in $\R^{n\times 1} = \R^n$, are assumed to be column vectors. For a vector $\bfx\in \mathbb{R}^d$ we use $\|\bfx\|$ to denote the $L_2$-norm of $\bfx$. A vector $\bfr\in \mathbb{R}^n$ is called stochastic if $r_i\geq 0$ and $\sum_{i=1}^n r_i=1$. Similarly, a non-negative matrix $A\in \mathbb{R}^{n\times d}$ is called (row) stochastic if $\sum_{j=1}^d A_{ij}=1$ for every $i\in [n]$. 
For a matrix $A\in \mathbb{R}^{n\times d}$, we denote its $i$-th row and $j$-th column by $A_i$ and $A^j$, respectively.  For an $n\times d$  matrix $A$ and a strictly positive stochastic vector $\bfr\in \mathbb{R}^n$, we define the $\bfr$-norm of $A$  by $\lnr A \rnr_{\bfr}^2 = \sum_{i=1}^n r_i \lnr A_i\rnr^2$. It can be verified that $\Nr{\cdot}$ is a norm on the space of $n\times d$ matrices. We denote the Frobenius norm of $A$ by $\lnr A\rnr_F$, where $\lnr A\rnr_F^2 = \sum_{i=1}^n \sum_{j=1}^d |A_{ij}|^2$. Moreover,  $A \geq B$ indicates that all the entries of $A-B$ are non-negative. 

\section{Problem Setup and Main Result}\label{sec:prob-def-main-res}
In this section, we discuss the problem formulation and the main result of this work. 
\subsection{Problem Setup}
Consider a set of $n\geq 2$ agents that are connected through a time-varying network. Each agent $i\in[n]$ has access to a local cost function $f_i:\mathbb{R}^{1\times d}\rightarrow \mathbb{R}$. The goal of this paper is to minimize the function $f(\bfx):= \sum_{i=1}^n r_i f_i(\bfx)$, or equivalently  solve 
\begin{align}\label{eq:ERM3}
     \min_{\bfx_1,\ldots,\bfx_n\in \mathbb{R}^d} \sum_{i=1}^{n}r_i f_i(\bfx_i) \quad \textrm{subject to}\quad \bfx_1=\bfx_2=\cdots=\bfx_n,
	\end{align}
	where vector $\bfr=(r_1,r_2,\dots, r_n)$ is a stochastic vector, i.e., $r_i\geq 0$ and $\sum_{i=1}^n r_i=1$. 
	
   We represent the time-varying topology at time $t\geq 1$ by the directed graph $\clg(t)=([n],\cle(t))$, where the vertex set $[n]$ represents the set of agents and the edge set ${\cle(t)\subseteq\{(i,j):i,j\in [n]\}}$ represents the set of links at time $t$. At each time $t$, agent $i$ can only send messages to its (out-) neighbors in $\cle(t)$, i.e., all $j\in [n]$ such that $(i,j)\in \cle(t)$. In order to achieve a consensus, the sequence $\{\clg(t)\}$ should satisfy some desirable long-term connectivity properties which will be discussed in Assumption~\ref{asm:W}.

   To present our algorithm for solving \eqref{eq:ERM3} collaboratively, let us first discuss the general framework for lossy/noisy sharing of information that is considered in this work. We assume that each agent maintains the sate $\bfx_i(t)\in \R^d$, which is an estimate of the optimizer of~\eqref{eq:ERM3}, and has access to its local cost function's gradient information. Moreover, it has access to an imperfect weighted average of its  neighbors states at time $t$, denoted by $\xh_i(t)$. More precisely, agent $i$ has access to $\xh_i(t)=\sum_{j=1}^nW_{ij}(t)\bfx_j(t)+\nt_i(t)$ where ${W(t)=[W_{ij}(t)]}$ is a row-stochastic matrix that is consistent with the underlying network $\clg(t)$ and $\nt_i(t)$ is a random noise vector in $\R^d$. By consistency of $W(t)$ and $\clg(t)$ we mean that that $W_{ij}(t)>0$ if only if $(j,i)\in\clg(t)$. Although it appears to be a simplistic model, later in Section~\ref{sec:x-h-implic} we provide some practical implications of such a structural assumption on $\xh_i(t)$ and its generality to contain many models of lossy and noisy sharing of information.

   Now, we are ready to present our \textit{Diminishing Mixing} (\texttt{DIMIX}) algorithm. In this algorithm, each agent $i$ updates its current estimate by computing a diminishing convex combination of its own state and received noisy average estimate $\xh_i(t)$, moving along its local gradient. More formally, the update rule is given as the following
   \begin{align}\label{eq:upd_nd}
    \bfx_i(t+1)\!=\! (1-\beta(t))\bfx_i(t)\!+\! \beta(t)\xh_i(t)\!-\! \alpha(t) \beta(t)\nabla f_i(\bfx_i(t)),
   \end{align}
   where $\alpha(t) = \frac{\azr}{t^{\nu}}$, $\beta(t) = \frac{\bzr}{t^\mu}$, and $\mu, \nu \in (0,1)$ are the diminishing step-sizes of the algorithm. The description of \texttt{DIMIX} is summarized in Algorithm~\ref{alg:1}. A similar algorithm is independently proposed and discussed for a particular subsetting (i.e., specific lossy sharing mechanism, weight vector $\bfr$, and a specific choice of $\nu,\mu>0$) of our framework and for time-invariant networks in \cite{vasconcelos2021improved}.

   For simplicity of notation, let 
     \begin{align}
        X(t) :=\begin{bmatrix} 
        \bfx_1(t) \\ \vdots \\ \bfx_n(t)
        \end{bmatrix},\quad 
        E(t) := \begin{bmatrix}
        \nt_1(t) \\ \vdots\\ \nt_n(t)\end{bmatrix},\quad 
        \nabla f(X(t)) := \begin{bmatrix} \nabla f_1(\bfx_1(t))\\ \vdots \\ \nabla f_n(\bfx_n(t))\end{bmatrix}.
        \label{eq:def:grad-f-X}
    \end{align}
    Using this notation, we can rewrite the update algorithm~\eqref{eq:upd_nd} in the  matrix format as
    \begin{align}
        X(t+1) & = ((1-\et)I+\et W(t))X(t)+\et\ER{t} - \at\et\nabla f(X(t)). 
        \label{eq:dyn-col-mat}
    \end{align}
    
\begin{algorithm}[t!]
    \caption{\texttt{DIMIX} at agent $i$}\label{alg:1}
    \begin{algorithmic}[1]
    \Require Stochastic matrix sequence $\{W(t)\}$, Iteration $T$
    \State Set $\bfx_{i}(1)=0$. 
    \For{$t=1,\ldots,T-1$}
    \State Compute the local gradient $\nabla f_i(\bfx_{i}(t))$.
    \State Obtain noisy average neighbors' estimate $\xh_i(t)$.
    \State Set: 
      $\bfx_{i}(t+1) = (1-\beta(t))\bfx_i(t)+ \beta(t)\xh_i(t)- \alpha(t) \beta(t)\nabla f_i(\bfx_i(t))$.
\EndFor
\end{algorithmic}
\end{algorithm}
 
\subsection{Assumptions} 
Here, we discuss the assumptions that we will use in the subsequent discussions in this work, namely, assumptions on the agent~$i$'s neighbor average state estimate $\xh_i(t)$, the stochastic weight matrix $\{W(t\})$, and local objective functions $f_i$.
 \begin{asm}[\textbf{Noise Assumption}]\label{assum:neighbor}
    We assume that the noise sequence $\{\nt_i(t)\}$ satisfies         \begin{align*}
            \EE{\nt_i(t)\mid \F_t}&=0\mbox{, and}\cr 
           \EE{\|\nt_i(t)\|^2 \mid \F_t}&\leq \sm,
        \end{align*}
        for some $\sm>0$, all $i\in[n]$, and all $t\geq 1$. Here, $\{\F_t\}$ is the natural filtration for the process $\{X(t)\}$. 
 \end{asm}
As mentioned before, to provide guarantees on the working of our algorithm,  certain connectivity assumptions need to be satisfied among the agents over time. 
\begin{asm}[\textbf{Connectivity Assumption}]\label{asm:W}
 We assume that the weight matrix sequence $\Wt$  satisfies the following properties.
 \begin{enumerate}[(a)]
     \item \textit{Stochastic with Common Stationary Distribution}: $W(t)$ is non-negative and  ${W(t)\ones=\ones}$ and ${\bfr^T W(t)=\bfr^T}$ for all $t\geq 1$, where $\ones\in \R^n$ is the all-one vector, and $\bfr>0$ is a given stochastic weight vector. 
     \item \textit{Bounded Nonzero Elements}: There exists some $\minW>0$ such that if for some $i,j\in [n]$ and $t\geq 1$ we have  $W_{ij}(t)>0$, then $W_{ij}(t)\geq \minW$.
     \item \textit{$B$-Connected}: For a fixed integer $B\geq 1$, the graph $\lp [n],\bigcup_{k=t+1}^{t+B}\cle(k)\rp$ is strongly connected for all $t\geq 1$, where $\cle(k)=\{(j,i)\mid W_{ij}(k)>0\}$. 
     
 \end{enumerate}
 \end{asm}
 
 For the local objective functions, we make the following assumption. 
 
\begin{asm}[\textbf{Function Assumptions}]\label{asm:f}
We assume the following properties on the function $f_i$ for all $i$
\begin{enumerate}[(a)]
     \item	The function $f_i$ is $\li$-smooth, i.e., for any $\bfx,\bfy\in\mathbb{R}^d$ we have that $\|\nabla f_i(\bfx)-\nabla f_i(\bfy)\|\leq \li\|\bfx-\bfy\|$.
    
    \item The function $f_i$ is $\m$-strongly convex, i.e., for any $\bfx,\bfy\in \R^d$, we have ${\la \nabla f_i(\bfx)\hspace{-1pt}-\hspace{-1pt}\nabla f_i(\bfy) , \bfx \hspace{-1pt}-\hspace{-1pt} \bfy \ra \geq \m \|\bfx\hspace{-1pt}-\hspace{-1pt}\bfy \|^2}$.
 \end{enumerate}
\end{asm}

\begin{remark}\label{rem:asm:f}
Since $\bfr$ is a stochastic vector, the properties of $f_i$s in Assumption~\ref{asm:f} can be immediately translated to similar properties for $f(\bfx)=\sum_{i=1}^n r_i f_i(\bfx)$. Thus, the function  $f$ is also $L$-smooth and $\m$-strongly convex.  
\end{remark}

\subsection{Main Result and Discussion}
The main result of this paper is the following theorem.
\begin{theorem}\label{thm:strng_cvx}
Assume the conditions in Assumptions~\ref{assum:neighbor}--\ref{asm:f} are satisfied and the step-sizes are set to  $\alpha(t) = \frac{\azr}{t^{\nu}}$ and $\beta(t) = \frac{\bzr}{t^\mu}$ for $\mu,\nu \in (0,1)$. Then, if $\mu+\nu<1$, the dynamics generated by Algorithm~\ref{alg:1} satisfy
\begin{align}\label{eq:main-1}
    \EE{\Nr{X(T)-\ones\bfx^\s}^2} &\leq
    \cOne T^{-\min(\mu,2\nu)}+ \cTwo T^{-\min(\mu-\nu,2\nu)},
\end{align}
for any iteration $T\geq \Ts:=\max(\Tzero,\Tfirst, \TStrt,\Ta,\Tb, \Tc, \Td, \Te)$, where $\Tzero$,
$\Tfirst$, $\TStrt$, $\Ta$, $\Tb$, $\Tc$, $\Td$ and $\Te$  are given in~\eqref{eq:def:T0}, ~\eqref{eq:def:Tfirst}, ~\eqref{eq:tau_0},~\eqref{eq:Ta},~\eqref{eq:Tb},~\eqref{eq:Tc},~\eqref{eq:Td}, and~\eqref{eq:Te},    respectively,  
and ${\bfx^\s :=\arg\min f(\bfx)}$.  Furthermore, under the same assumptions, when $\mu+\nu=1$, 
the dynamics generated by Algorithm~\ref{alg:1} satisfy~\eqref{eq:main-1}, 
for any iteration $T\geq \Ts$, provided that $\frac{\rho \li}{\rho +\li} \azr \bzr \geq 8 \min(\mu-\nu, 2\nu)$.
\end{theorem}

We refer to  Section~\ref{sec:proof} for the  proof of Theorem~\ref{thm:strng_cvx}. 

\begin{remark}
Theorem~\ref{thm:strng_cvx} guarantees the exact convergence (in $L_2$ sense) of each local state to the global optimal with diminishing step-size even though the noises induced by random quantizations and gradients are non-vanishing with iterations. In order to maximize the exponents in the upper bound~\eqref{eq:main-1}, 
it can be verified that the optimum choice is $(\mu,\nu) = (3/4,1/4)$. Replacing this in~\eqref{eq:main-1}, we conclude 
\begin{align*}
    &\EE{\Nr{X(T)-\ones\bfx^\s}^2} \leq  \xi T^{-1/2},
\end{align*}
for any {$T\geq \Ts$} and $\xi = \cOne+\cTwo$. Our algorithm and the main result are inspired by the fixed step-size variation of \eqref{eq:dyn-col-mat} that is proposed  in~\cite{reisizadeh2019exact} under the limited setting of time-invariant networks, uniform weights $\bfr$, and a particular choice of lossy sharing of information. In that setting, it is shown that for any \textit{given stopping time} $T\geq T_{\min}$ and any $\epsilon>0$, the constant step-sizes $\alpha_0,\beta_0>0$ can be set such that \[{\EE{\Nr{X(T)-\ones\bfx^\s}^2}\leq  c T^{-1/2+\epsilon,}}\]
where $c,T_{\min}$ are positive constants depending on the problem's parameters (note that this is established for a fixed $T$). However, $T_{\min} \rightarrow \infty$ as $\epsilon \rightarrow 0$~\cite{reisizadeh2019exact}. Here, we provide a rigorous convergence rate analysis which reduces to $\BO\lp T^{-1/2}\rp$ for \textit{every} iteration $T$. In Theorem~\ref{thm:strng_cvx}, for the case $\mu+\nu=1$, the minimum number of required iterations is finite.

We also note that in a recent independent work~\cite{vasconcelos2021improved}, for a specific quantizer and the specific choice $\bfr=\frac{1}{n}\ones$, the authors have shown the convergence rate of $\BO(T^{-1/2}(\ln T)^2)$ for strongly convex loss functions over \textit{fixed} underlying networks and a specific choice of lossy sharing of information. Note that the obtained rate, which is with respect to a weighted average of the previous iterates $\{X(t)\}_{t\leq T}$ instead of $X(T)$, is strictly \textit{slower} than $\BO(T^{-1/2})$. Furthermore, it is assumed that gradients are bounded in~\cite{vasconcelos2021improved}, while here we show that such a strong assumption is not needed and Lipschitz gradients result in expected bounded gradients for strongly convex functions.
\end{remark}

One of the key distinctions of our work from the prior works in this domain is the introduction of $\hat{x}_i(t)$ that satisfies some general structural properties without being tied to any specific application. Before discussing the technical details of the main result's proof, let us provide some general practical settings for which our structural assumptions on $\hat{x}_i(t)$ hold. 

\subsection{Examples for Stochastic Noisy State Estimation}\label{sec:x-h-implic}
The noisy estimation of the neighbors' state used in~\eqref{eq:upd_nd} may appear in various practical problem settings, due to the physical limitations of communication between the computing agents in the network. In the following, we describe a few scenarios, in which each neighbor can only receive a noisy estimate of the updated states in each iteration. 
\begin{example}
	\textbf{(Noisy Communication).}
	A practical scenario in which the noisy neighbor estimate model will be realized is due to the communication between the agents over a noise channel.  Consider a wireless medium, in which the computing nodes communicate to their neighbors over a Gaussian channel, i.e., when node $j$ sends its state $\bfx_j(t)$ to its neighbor $i$,  the signal received at node $i$ is $\bfx_j(t) + \bfz_{i,j}(t)$ where $\bfz_{i,j}(t)$ is a zero-mean Gaussian noise with variance $\sigma^2$, independent across $(i,j)$, and $t$. Then we have 
	\begin{align*}
	    \hat{\bfx}_i(t) & = \sum_{j=1}^n W_{ij}(t) \lp \bfx_j(t) + \bfz_{i,j}(t)\rp\cr
	    & = \sum_{j=1}^n W_{ij}(t)  \bfx_j(t) + \sum_{j=1}^n W_{ij}(t)  \bfz_{i,j}(t).  
	\end{align*}
	Thus, we have $\nt_i(t) \hspace{-1pt}=\hspace{-1pt} \sum_{j=1}^n W_{ij}(t)  \bfz_{i,j}(t)$, which implies $\EE{\nt_i(t)}=0$ and ${\EE{\|\nt_i(t) \|^2} \hspace{-1pt}= \sigma^2 \sum_{j=1}^n W_{ij}(t)^2 \hspace{-1pt}\leq \sigma^2}$. 
	Hence, the conditions of Assumption~\ref{assum:neighbor} are satisfied.

\end{example}
	
\begin{example}\label{exm:SQb}\textbf{(Unbiased Stochastic Quantizer).} 
    	In many applications, there is band-limited link between the agents. That is, the state $\bfx$ of a user $i$ needs to be quantized to a certain number bits, before transmission to its neighbors. The difference between the actual state and its quantized version can be modelled as the estimation noise. In particular, the stochastic quantizer with a number of  bits $b$ maps a vector $\bfx\in \mathbb{R}^d$ to a random vector  $Q_s^{S}(\bfx) \in \R^d$,  where its $j$-th entry is given by
    	\begin{align}
    	\left[Q^{S}_s(\bfx)\right]_j = \frac{1}{s} \|\bfx\|\cdot \sign{x_j}\cdot \zeta\lp \frac{|x_j|}{\|\bfx\|},s \rp ,\quad  j\in[d].
    	\end{align}
    	Here, for $t\in[0,1]$ we define $\zeta(t, s) =  \lfloor s t \rfloor +\mathbbm{1}[U<st -  \lfloor s t \rfloor]$, where $U$ is random variable with uniform distribution over $[0,1]$, and $\mathbbm{1}[\cdot]$ is the binary indicator function. Note that the random variable $U$s are chosen to evaluate the $\zeta(\cdot,\cdot)$ function are independent, across the coordinates, agents, and time slots. 
    	It is shown in~\cite{alistarh2017qsgd} that applying this quantizer on $\bfx \in \mathbb{R}^d$ with a bounded norm $\|\bfx\|\leq D$ satisfies ${\EE{Q^{S}_s(\bfx)}=\bfx}$ and ${\EE{\|Q^{S}_s(\bfx) - \bfx\|^2}\leq \min\left(\frac{\sqrt{d}}{s}, \frac{d}{s^{2}}\right)D}$.
    	Therefore, the neighbors estimate for node $i$ will be 
    	\begin{align*}
    	    \hat{\bfx}_i(t) & = \sum_{j=1}^n W_{ij}(t) Q^{S}_s(\bfx_j(t)) \cr 
    	    & = \sum_{j=1}^n W_{ij}(t) \bfx_j(t) + \sum_{j=1}^n W_{ij}(t) \lp Q^{S}_s(\bfx_j(t)) - \bfx_j(t)\rp, 
    	\end{align*}
    	where
    	$\nt_i(t) = \sum_{j=1}^n W_{ij}(t) \lp Q^{S}_s(\bfx_j(t)) - \bfx_j(t)\rp$ is the estimation noise. From the properties of $Q^{S}_s(\cdot)$ and independence of $U$s, we can verify that $\EE{\nt_i(t)|\F_t} =0$ and 
    	\begin{align*}
    	   \EE{\|\nt_i(t)\|^2|\F_t} & = \min\left(\frac{\sqrt{d}}{s}, \frac{d}{s^{2}}\right)D \sum_{j=1}^n W_{ij}(t)^2 \cr
    	   &\leq \min\left(\frac{\sqrt{d}}{s}, \frac{d}{s^{2}}\right)D.
    	\end{align*}
    
    	Therefore, the conditions of Assumption~\ref{assum:neighbor} are satisfied. Note that here $\nt_i(t)$ and $\nt_j(t)$ might be correlated, especially if nodes~$i$ and~$j$ have common neighbors. However, this is not in contradiction with the assumption. 
	\end{example}
	
\section{Experimental Results}\label{sec:Exp-result}
	Here, we provide experimental results supporting the effectiveness of the proposed
    algorithms. We use linear regression for the validation of our algorithm in the strongly convex setting. 
    
    \noindent\textbf{Data and Experimental Setup.} We consider a $25$-dimensional linear regression problem. We synthesize $N=100$ data points $\{\omega_1,\ldots,\omega_{100}\}$, where $\omega_i=(\bfu_i,v_i)$ with $v_i=\bfu_i^T\Tilde{\bfx}+\theta_i$. In order to generate the data, we uniformly and independently draw the entries of each $\bfu_i$ and each $\theta_i$ from $(0,1)$ and $(0, 0.1)$, respectively. Similarly, entries of $\Tilde{\bfx}$ are sampled uniformly and independently from $(0,0.8)$. We consider the following loss function
    \[f(\bfx): = \frac{1}{2N}\sum_{i=1}^{N}\|v_i-\bfu_i^T\bfx\|^2.\]
    We set $n=20$ as the number of worker nodes and distribute the data points across the nodes according to $r_i=p_i/\sum_{i=1}^{20} p_i$, where $p_i$ is drawn uniformly at random from the interval $(0.01,0.09)$. 
    We implement the unbiased stochastic quantizer with ${s=4}$ (Example~\ref{exm:SQb}). We utilize the following stochastic matrix sequences for the simulations. 
    
    \noindent\textbf{Experiments with Fixed Graph.}
     We consider a fixed undirected cyclic graph on $\clg^C = ([n], \mathcal{E})$, where ${\mathcal{E} = \{(\la i \ra,\la i+1\ra): i\in [n]\}}$, where $\la i \ra = (i \mod n)+1$.  
     The stochastic matrix sequence $W(t)=W$ for any $t$ is given by
    \begin{align}\label{W:fdx}
          W_{ij}=\!\! \begin{cases}
          \frac{r_{\la j \ra}}{2(r_{\la i \ra} +r_{\la j \ra})} & j \in\{\la i-1 \ra, \la i+1 \ra\} \\
          \frac{r_{\la i \ra}}{2(r_{\la i \ra}+r_{\la i+1 \ra})}\!+\! \frac{r_{\la i \ra}}{2(r_{\la i \ra}+r_{\la i-1 \ra})} & j = i \\
          0 & \textrm{otherwise}.
          \end{cases}
     \end{align}
     
    \noindent\textbf{Experiments with Time-varying Graph.}
    To evaluate Algorithm~\ref{alg:1} for a time-varying graph, we implement a class of cyclic gossip algorithms~\cite{boyd2006randomized,aysal2009broadcast,kar2012distributed}. This algorithm goes through cycles of the above fixed graph $\clg^C$ where in each iteration only one pair of neighbor nodes exchange information. 
    
    More precisely, we utilize a family of $n$ undirected graphs, each including a single edge from the cyclic graph above. We have $\mathcal{G}(t) = ([n], \mathcal{E}(t))$, where $\mathcal{E}(t) = \{({\la t \ra}, {\la t+1 \ra} )\}$. The  stochastic matrix sequence $\{W(t)\}$ corresponding to $\mathcal{E}(t)$ is given by
    \begin{align}\label{W:var}
        [W(t)]_{ij} =
        \begin{cases}
        \frac{r_{\la j\ra}}{ r_{{\la t\ra}} + r_{{\la t+1\ra}}} & i,j\in \{{\la t\ra}, {\la t+1\ra}\}\\
        1 & i=j \notin \{{\la t\ra}, {\la t+1\ra}\}\\
        0 & \textrm{otherwise}.
        \end{cases}
    \end{align}
    Note that these $n$ edges will be experienced in a periodic manner, and the properties of  Assumption~\ref{asm:W} holds for~\eqref{W:fdx} and~\eqref{W:var}. In particular, the sequence of stochastic matrices $\{W(t)\}$ is $n$-connected. 
    
    \begin{figure}[h]
    \centering
    \includegraphics[width=0.50\linewidth]{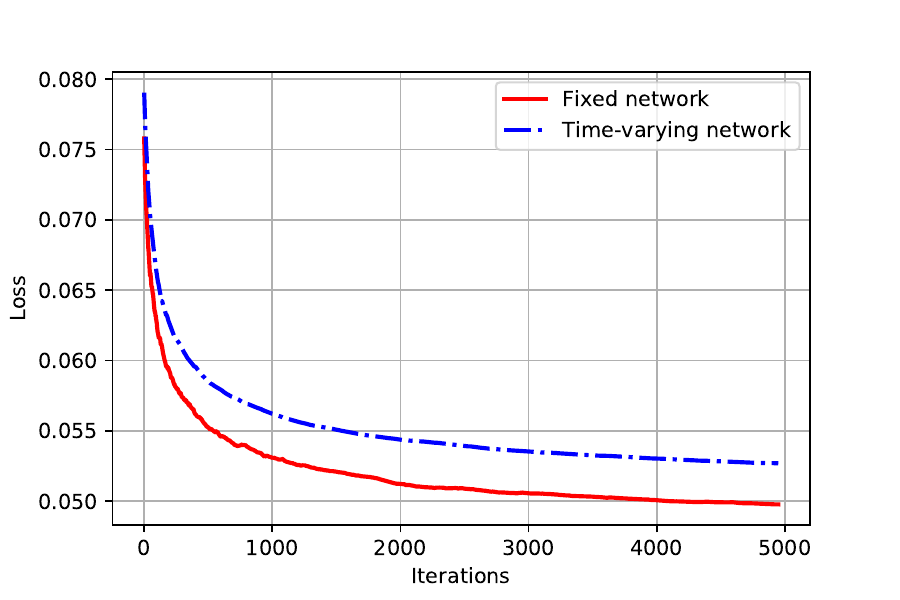}
    \caption{Training Loss vs. Iterations: Linear Regression on the synthetic data.}\label{fig:fig-loss}
    \includegraphics[width=0.50\linewidth]{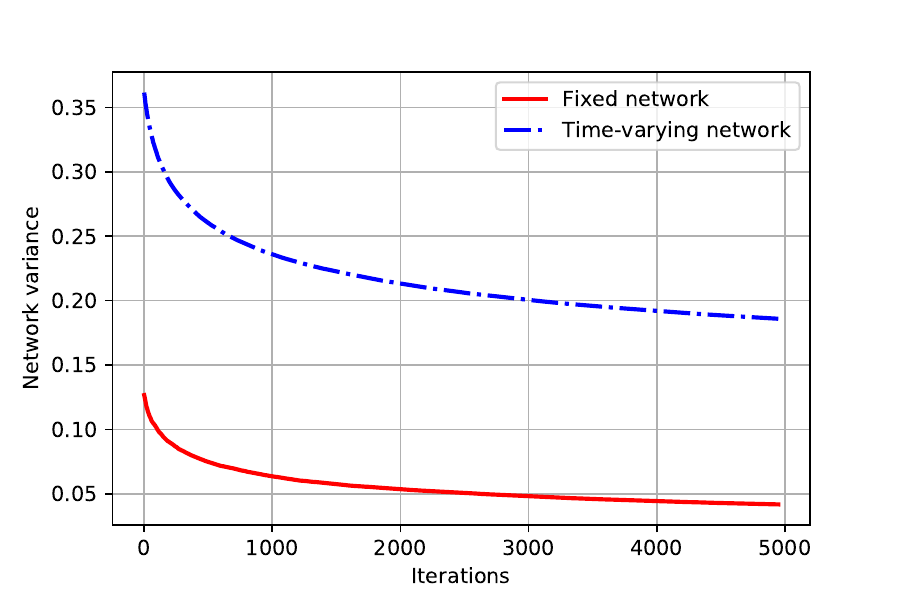}
    \caption{Network variance vs. Iterations: Linear Regression on the synthetic data.}\label{fig:fig-var}
    \end{figure}
    
    Figure~\ref{fig:fig-loss} and~\ref{fig:fig-var} demonstrate the training time vs. the objective loss function and the training time vs. the deviation of each node's state from the average state of a linear regression, respectively for the synthetic dataset over the networks introduced above. Here, `Fixed network' refers to the full cycle with stochastic matrix in~\eqref{W:fdx} and `Time-varying network' refers to the network with stochastic matrix sequence in~\eqref{W:var}. The experiments are performed for $T=5000$ iterations of SGD. The parameters of the dynamics in~\eqref{eq:upd_nd} are finely tuned to $(\azr,\nu^\s) = (0.1,1/4)$ and $(\bzr,\mu^\s) = (0.7,3/4)$ for the both discussed fixed and time-varying graphs. It can be verified that a faster mixing of information over the network in the same architecture leads to a faster convergence.
    
   \section{Auxiliary Lemmas}\label{sec:aux-lemmas}
   In this section, we present auxiliary lemmas which play crucial roles in the proof of the  main result, namely, Theorem~\ref{thm:strng_cvx} in Section~\ref{sec:proof}. The proofs of Lemmas~\ref{lemma:transition2}--\ref{lm:sum_dl_gm} are provided in Appendix. We refer to~\cite{nesterov1998introductory} for the proof of Lemma~\ref{lm:sc_bnd}.
   
\begin{lemma}\label{lemma:transition2}
    Let $\Wt$ satisfy the connectivity Assumption~\ref{asm:W} with parameters $(B,\minW)$,  and let  $\{A(t)\}$ be given by ${A(t)=(1-\beta(t))I+\beta(t)W(t)}$ where $\beta(t)\in (0,1]$ for all $t$, and $\{\beta(t)\}$ is a non-increasing sequence. Then, for any matrix $U\in \mathbb{R}^{n\times d}$, and all $t>s\geq 1$, we have 
    \begin{align*}
    &\Nr{\lp A(t-1)A(t-2)\cdots A(s+1)-\ones \bfr^T\rp U }^2 \leq \kap \prod_{k=s+1}^{t-1}(1-\lambda\beta(k))\Nr{U}^2,
    \end{align*}
    where  $\lambda:=\frac{\minW \rmin}{2Bn^2}$,  $\kap := \lp 1-B\lambda \bzr\rp^{-1}$ and $\bzr=\beta(1)$. 
\end{lemma}

\begin{lemma}\label{lemma:matrixprop}
    For an $n\times m$ matrix $A$ and $m\times q$ matrix $B$, we have
    \begin{align*}
         \lnr AB \rnr_{\bfr}	   \leq     \lnr A \rnr_{\bfr} \lnr B \rnr_{F}.  
    \end{align*}
\end{lemma}
	
\begin{lemma}\label{lm:nor_squr}
   For any pair of vectors $\bfu$, $\bfv$, and any scalar $\theta>0$, we have 
   \begin{align*}
        \|\bfu+\bfv\|^2 &\leq (1+\theta)  \|\bfu\|^2 + \lp 1+ \frac{1}{\theta} \rp \|\bfv\|^2.
   \end{align*}
   Similarly, for matrices $U$ and $V$ and any scalar $\theta>0$, we get 
   \begin{align*}
       \Nr{U+V}^2 \leq (1+\theta)\Nr{U}^2 + \lp 1 + \frac{1}{\theta} \rp\Nr{V}.
   \end{align*}
\end{lemma}
 \begin{lemma}\label{lm:delta}
        For any $0\leq \dl < 1$ and $0<a<1$  we have
        \begin{align*}
             \prod_{k=s}^{t-1} \lp 1-\frac{a}{k^{\dl}} \rp 
         &
         \leq \exp\lp -\frac{a}{1-\dl} \lp t^{1-\dl} - s^{1-\dl}\rp \rp.
        \end{align*}
        For $\dl = 1$ and $0\leq a< 1$  we have 
        \begin{align*}
             \prod_{k=s}^{t-1} \lp 1-\frac{a}{k} \rp 
         &
         \leq \lp\frac{t}{s}\rp^{-a}.
        \end{align*}
    \end{lemma}
     	\begin{lemma}\label{lemma:Psum}
    	Let $\{\beta(t)\}$ be a sequence in $\R$ and $\lambda$ be a non-zero scalar. Then the following identities hold for all $ t\geq 1$
    	 \begin{align}\label{eq:esumprod}
    	    &\sum_{s=1}^{t-1} \beta(s)\prod_{k=s+1}^{t-1}(1-\lambda \beta(k))=\frac{1}{\lambda}-\frac{1}{\lambda}\prod_{k=1}^{t-1}(1-\lambda \beta(k)),
         \end{align}
    As a result, for any sequence $\{\beta(t)\}$ in $[0,1]$ and $\lambda>0$, we get
    \begin{align*}
      \sum_{s=1}^{t-1} \beta(s)\prod_{k=s+1}^{t-1}(1-\lambda \beta(k)) \leq \frac{1}{\lambda}. 
    \end{align*}
    \end{lemma}
\begin{lemma}\label{lm:sum_dl_gm}
        For any $0\leq \dl < \min(1,\sigma)$, $0< a\leq 1$, and every
        ${t> \tLM := (\frac{2(\sigma -\dl )}{a})^{\frac{1}{1-\dl}}}$, we have
        \begin{align*}
            \sum_{s=1}^{t-1}  \lb \frac{1}{s^{\sigma}} \prod_{k=s+1}^{t-1} \lp 1-\frac{a}{k^{\dl}} \rp \rb 
         &
         \leq A(a,\sigma,\dl)t^{-(\sigma-\dl)},
        \end{align*}
        where $A(a,\sigma,\dl)$ is given by
        \begin{align}\label{eq:def:A}
        A(a,\sigma,\dl):=\begin{cases} 
        2^\sigma \max\left\{1+\frac{2}{a},1+\frac{1}{\sigma-1} \left(\frac{2(\sigma-\delta)}{a}\right)^{\frac{\sigma-\delta}{1-\delta}}\right\} &\mbox{if } \sigma>1, \\
       2^\sigma \max\lc 1+ \frac{2}{a}, 1+ \frac{2}{a}\ln\lp\frac{2(1-\dl)}{a}\rp\rc & \mbox{if } \sigma=1,\\
         2^\sigma \max\lc 1+ \frac{2}{a}, 1+ \frac{2(\sigma-\delta)}{a(1-\sigma)}\rc &\mbox{if } 0<\sigma<1.
       \end{cases}
        \end{align}
        Moreover, for $\dl=1$ and $a-\sigma +1 \neq 0$, we have
        \begin{align*}
            \sum_{s=1}^{t-1}  \lb \frac{1}{s^{\sigma}} \prod_{k=s+1}^{t-1} \lp 1-\frac{a}{k} \rp \rb \leq A(a,\sigma,1) t^{-\min(\sigma-1,a)},
        \end{align*}
        where $A(a,\sigma,1) = 2^\sigma \lp 1+ \frac{1}{|a-\sigma +1|}\rp $.
    \end{lemma}
    
    \begin{lemma}\label{lm:t->t+1}
    For non-negative numbers $a,b,c>0$, if $b\neq 1$ we have 
    \begin{align}\label{eq:lm:t->t+1}
        (1-a t^{-b}) t^{-c} \leq (t+1)^{-c}
    \end{align}
    for every $t\geq t_0 :=\lp \frac{c}{a}\rp^{\frac{1}{1-b}}$. Moreover, if $b=1$, the inequality in~\eqref{eq:lm:t->t+1} holds for every $t\geq 1$ provided that $a\geq c$.
    \end{lemma}
    
       \begin{lemma}\label{lm:sc_bnd}\cite[Theorem 2.1.11]{nesterov1998introductory}
      Suppose that $\nabla f$ is Lipschitz continuous with constant $L$ and $f$ is strongly convex with modulus $\m$. Then, we have
     \begin{align*}
      \la \bfx-\bfx^\s, \nabla f(\bfx)  \ra \geq c_1\|\nabla f(\bfx)\|^2 +c_2 \|\bfx-\bfx^\s\|^2,
     \end{align*}
     where $c_1=\frac{1}{\m+L}$, $c_2=\frac{\m L}{\m + L}$, and $\nabla f(\bfx^\s)=0$.
    \end{lemma}

	\section{Proof of Theorem~\ref{thm:strng_cvx}}
	\label{sec:proof}
	In this section, we provide the proof of the main result, namely, Theorem~\ref{thm:strng_cvx}. The main steps of the proof are twofold: We first bound the deviation of the agents' states from their average, and then  analyze the distance of the average state from the global optimal point. These together lead to the proof of the theorem. 
	\subsection{State Deviation from the Average State}
	\label{sec:node-to-avg}
	In this part, we study $\C(t):=\EE{\sum_{i=1}^n r_i \lnr\bfx_i(t)- \mx(t)\rnr^2}$, where $\mx(t) :=\sum_{i=1}^n r_i \bfx_i(t)$ is the avarage of the states at time $t$. 
	Note that the dynamics in~\eqref{eq:dyn-col-mat} can be viewed as the linear time-varying system  
	\begin{align}\label{eq:controlled}
	    X(t+1)=A(t)X(t)+U(t),
	\end{align}
    with
    \begin{align*}
    A(t)&=((1-\et)I+\et W(t)),\cr
    U(t)&=\et\ER{t} - \at\et\nabla f(X(t)). 
    \end{align*}
     Therefore, we have
	\begin{align}
	    X(t)=\sum_{s=1}^{t-1}\Phi(t:s)U(s) + \Phi(t:0)X(1),
	    \label{eq:Z(t)}
	\end{align}
	where $\Phi(t:s)=A(t-1)\cdots A(s+1)$ with $\Phi(t:t-1)=I$ is the transition matrix of the linear system \eqref{eq:controlled}. 
	For the notational simplicity, let us define 
	\begin{align*}
	    P(t:s)&:=\beta(s)(\Phi(t:s)-\ones\bfr^T) \cr
	    & = \beta(s)\lp A(t-1)\cdots A(s+1)-\ones\bfr^T\rp.
	\end{align*}
     As a result of Lemma~\ref{lemma:transition2}, we have
	     \[\Nr{ P(t:s)U}\leq \pi(t:s)\Nr{U},\]
	     where $\pi(t:s)$ is defined by 
	     \begin{align}\label{eqn:pidef}
	         \pi(t:s):=\beta(s)\kap^{\frac{1}{2}}\prod_{k=s+1}^{t-1}(1-\lambda\beta(k))^{\frac{1}{2}}.
	     \end{align}
	Assuming $X(1)={\mathbf{0}}$, the dynamic in~\eqref{eq:Z(t)} reduces to 
	\begin{align}
	    X(t)=\sum_{s=1}^{t-1}\Phi(t:s)U(s).
	    \label{eq:Z2(t)}
	\end{align}
	Moreover,  multiplying both sides of \eqref{eq:Z2(t)} from the left by $\bfr^T$ and using the fact that $\bfr^TA(t)=\bfr^T$, we get 
	\begin{align}
	\mx(t):=\bfr^T X(t) = \sum_{s=1}^{t-1}\bfr^T \Phi(t:s)U(s) = \sum_{s=1}^{t-1}\bfr^T U(s).
	\label{eq:Z:mean}
	\end{align}
	Then, subtracting~\eqref{eq:Z:mean} from~\eqref{eq:Z2(t)}, and plugging the definition of $U(s)$ we have
	\begin{align*}
	     X(t)-\ones \bar{\bfx}(t)
	    &= \sum_{s=1}^{t-1}\Phi(t:s)U(s) - \sum_{s=1}^{t-1} \ones \bfr^T U(s)\cr
	    &=\sum_{s=1}^{t-1}(\Phi(t:s)-\ones\bfr^T)U(s)\cr
	    &=\sum_{s=1}^{t-1}\beta(s)(\Phi(t:s)-\ones\bfr^T)\bigg[ \ER{s} - \alpha(s)\nabla f(X(s))\bigg]\cr
	    &=\sum_{s=1}^{t-1} P(t:s)\ER{s} -\sum_{s=1}^{t-1}\alpha(s)P(t:s)\nabla f(X(s)).
	\end{align*}

	Using Lemma~\ref{lm:nor_squr}  with $\theta=1$, we get
    \begin{align}
	    \Nr{X(t)-\ones\bar{\bfx}(t)}^2	
	    & \leq 2\Nr{\sum_{s=1}^{t-1} P(t:s)\ER{s}}^2 \!\!+\! 2\Nr{\sum_{s=1}^{t-1}\alpha(s)  P(t:s)\nabla f(X(s))}^2.\nonumber\\
	    & =2\sum_{s=1}^{t-1}  \Nr{ P(t:s)\ER{s}}^2 
	     + 2\sum_{s\neq q} \la P(t\hspace{-1pt}:\hspace{-1pt}s)\ER{s}, P(t\hspace{-1pt}:\hspace{-1pt}q)\ER{q}\ra \nonumber\\
	     &\phantom{=} + 2\Nr{ \sum_{s=1}^{t-1} \alpha(s)   P(t\hspace{-1pt}:\hspace{-1pt}s)\nabla f(X(s))}^2.\label{eq:Mtmain}
	\end{align}
	Using facts that $\ER{s}$ is measurable with respect to $\F_q$ for $q>s$ and $\EE{\ER{q}\md \F_q} = 0$,  we have
	\begin{align}\label{eq:inner-zero}
        \EE{\la P(t:s) \ER{s}, P(t:q) \ER{q}\ra}
        & =	
        \EE{\EE{\la P(t:s) \ER{s}, P(t:q) \ER{q}\ra \md \F_q}}\nonumber\\
        & =\EE{\la P(t:s) \ER{s}, P(t:q) \EE{\ER{q}\md \F_q}\ra} =0. 
    \end{align}
    Using a similar argument for $q<s$ and conditioning on $\F_s$, we conclude that~\eqref{eq:inner-zero} holds for all $q\not=s$. Recall that we defined  $\C(t)=\EE{\Nr{X(t)-\ones \mx(t)}^2}$. Taking  expectations of both sides of \eqref{eq:Mtmain} and using the identity in~\eqref{eq:inner-zero},  we get
	\begin{align}
	    C(t)\!=\!\EE{\Nr{X(t)-\ones \bar{\bfx}(t)}^2}
	    \leq 2\sum_{s=1}^{t-1} \EE{ \Nr{ P(t:s)\ER{s}}^2}  +2\EE{ \Nr{\sum_{s=1}^{t-1} \alpha(s)  P(t:s)\nabla f(X(s))}^2}.
	    \label{eq:Mtmain:exp} 
	\end{align}
    
    We continue with bounding the first term in~\eqref{eq:Mtmain:exp}. 
   From Assumption~\ref{assum:neighbor}, we have
   \begin{align*}
      \EE{\Nr{E(s)}^2}  = \EE{\EE{\Nr{E(s)}^2|\F_s}}
       = \EE{ \sum_{i=1}^n r_i \EE{\|\bfe_i(s)\|^2 |\F_s}}\leq \EE{\sum_{i=1}^n r_i \gamma}= \gamma.
   \end{align*}
   This together with Lemma~\ref{lemma:transition2} lead to
    \begin{align}\label{eq:X-mx-term} 
    \sum_{s=1}^{t-1} \EE{ \Nr{ P(t:s)\ER{s} }^2}
   & \leq  \sum_{s=1}^{t-1} \lb \beta^2(s)\kap \prod_{k=s+1}^{t-1} (1-\lambda \beta(k)) \EE{ \Nr{ \ER{s}}^2} \rb \cr   
   & \leq \sm \kap \sum_{s=1}^{t-1} \lb \beta^2(s) \prod_{k=s+1}^{t-1} (1-\lambda \beta(k)) \rb\nonumber\\
   & = \sm \kap\sum_{s=1}^{t-1} \lb \frac{\bzr^2}{s^{2\mu}}\prod_{k=s+1}^{t-1} \lp 1-\lambda \frac{\bzr}{k^\mu} \rp \rb \leq \cvxTwo t^{-\mu},
   \end{align}
   where the last inequality holds for $t\geq \Tzero$, where 
   \begin{align}\label{eq:def:T0}
       \Tzero:= \left\lceil \lp 2\mu/\lambda \bzr \rp^{\frac{1}{1-\mu}} \right\rceil,
   \end{align}
   and follows from   Lemma~\ref{lm:sum_dl_gm}, with parameters $(\sigma, \delta, \tLM)=(2\mu, \mu, \Tzero)$. Moreover,  we have  ${\cvxTwo:= \sm \kap\bzr^2 A(\lambda \bzr, 2\mu, \mu)}$.

   Next, we bound the second term in~\eqref{eq:Mtmain:exp}. Note that  $\Nr{\cdot}$ is a norm. Hence,  using the triangle inequality we have 
\begin{align}\label{eq:Egradient}
	\!\!\!\!\!\!\EE{ \Nr{\sum_{s=1}^{t-1} \alpha(s)  P(t:s)\nabla f(X(s))}^2}
	&\!\leq \EE{ \lp\sum_{s=1}^{t-1} \Nr{\alpha(s)  P(t:s)\nabla f(X(s))}\rp^2}\cr 
	&\!=\!\!\!\!\!\sum_{1\leq s,q \leq t-1} \!\!\!\!\!\mathbb{E}\Big[ \alpha(s) \Nr{ P(t:s)\nabla f(X(s))} \alpha(q) \Nr{ P(t:q)\nabla f(X(q))} \!\Big].\ \ \ 
\end{align}
	  Using Lemma~\ref{lemma:transition2} and the fact that $2ab\leq a^2+b^2$, we can upper-bound this expression as

	 \begin{align}\label{eq:term4-2}
	 \sum_{1\leq s, q \leq t-1} \mathbb{E} \Big[ \alpha(s) &\Nr{ P(t:s)\nabla  f(X(s))} \cdot \alpha(q) \Nr{ P(t:q)\nabla f(X(q))}\Big]\cr
	  & \leq \sum_{1\leq s, q\leq t-1} 
	  \mathbb{E} \Big[ \alpha(s) \pi(t:s) \Nr{\nabla f(X(s))}   \cdot \alpha(q) \pi(t:q) \Nr{\nabla f(X(q))}\Big]  \cr
	  &=  \sum_{1\leq s, q\leq t-1}   \pi(t:s) \pi(t:q)\mathbb{E} \Big[ \alpha(s)\Nr{\nabla f(X(s))}   \cdot \alpha(q) \Nr{\nabla f(X(q))}\Big] \cr
	  &\leq \frac{1}{2} \sum_{1\leq s, q\leq t-1}  \pi(t:s) \pi(t:q) \mathbb{E}\Big[ \alpha^2(s)\Nr{\nabla f(X(s))}^2  + \alpha^2(q) \Nr{\nabla f(X(q))}^2\Big] \cr
	  &= \sum_{1\leq s, q\leq t-1}  \pi(t:s) \pi(t:q) \EE{\alpha^2(s) \Nr{\nabla f(X(s))}^2}\nonumber \\
	  &= \lp \sum_{q=1}^{t-1} \pi(t:q)\rp \cdot
	\lp \sum_{s=1}^{t-1}  \alpha^2(s) \pi(t:s) \EE{\Nr{\nabla f(X(s))}^2}\rp,
	\end{align}
	where $\pi(t:s)$ is given by \eqref{eqn:pidef}. Then, using the fact that ${\sqrt{1-x} \leq 1-x/2}$ and Lemma~\ref{lemma:Psum}, we get  
	\begin{align}
    	\sum_{q=1}^{t-1} \pi(t:q) &= \sum_{q=1}^{t-1}  \lb \beta(q) \kap^{\frac{1}{2}}     \prod_{k=q+1}^{t-1}  (1-\lambda \beta(k))^{\frac{1}{2}} \rb\cr 
    	& \leq \sum_{q=1}^{t-1} \beta(q)\kap^{\frac{1}{2}}   \prod_{k=q+1}^{t-1} \lp 1-\frac{\lambda}{2} \beta(k)\rp \leq \frac{2}{\lambda}\kap^{\frac{1}{2}}  .
    	 \label{eq:sum-Aq}
	\end{align}
	Moreover, we can write
	 \begin{align}\label{eq:normgrad}
     \mathbb{E}\!\left[\hspace{-1pt}\Nr{\nabla f(X(s))}^2\right] \!& = \EE{\Nr{\nabla f(X(s)) - \nabla f(\ones\mx(s)) + \nabla f(\ones\mx(s)) - \nabla f(\ones\bfx^\star) + \nabla f(\ones\bfx^\star)  }^2}\nonumber\\
     & \stackrel{\rm{(a)}}{\leq} \EE{3\Nr{\nabla f(X(s)) - \nabla f(\ones\mx(s))}^2 + 3 \Nr{\nabla f(\ones\mx(s)) - \nabla f(\ones\bfx^\star)}^2 + 3 \Nr{\nabla f(\ones\bfx^\star)}^2} \nonumber\\
     &= 
     3 \EE{\Nr{\nabla f(X(s)) - \nabla f(\ones\mx(s))}^2}  +3\EE{  \Nr{\nabla f(\ones\mx(s)) - \nabla f(\ones\bfx^\star)}^2} + 3 \EE{ \Nr{\nabla f(\ones\bfx^\star)}^2} \nonumber\\
     &= \!3\mathbb{E}\!\left[ \sum_{i=1}^n \!r_i\! \lnr\nabla f_i(\bfx_i(t)) \!-\! \nabla f_i(\mx(s))\rnr^2\right]
     \!\!+\! 3\mathbb{E}\!\left[ \sum_{i=1}^n \!r_i\! \lnr\nabla f_i(\mx(s)) \!-\! \nabla f_i(\bfx^\star)\rnr^2
     \right]
     \!\!+\! 3  \Nr{\nabla f(\ones\bfx^\star)}^2
     \nonumber\\
     & \stackrel{\rm{(b)}}{\leq} 
     3 \EE{\sum_{i=1}^n r_i L^2 \lnr\bfx_i(t) - \mx(s)\rnr^2}
     + 3 \EE{\sum_{i=1}^n r_i L^2 \lnr \mx(s) - \bfx^\star\rnr^2} +  3  \Nr{\nabla f(\ones\bfx^\star)}^2
     \nonumber\\
     &= 3 \li^2 \EE{ \Nr{ X(t) - \ones\mx(s)}^2}
     + 3 \li^2 \EE{\lnr \mx(s) - \bfx^\star\rnr^2} +  3  \Nr{\nabla f(\ones\bfx^\star)}^2
     \nonumber\\
     &=3 \li^2 \lp \C(s) + \Q(s) \rp  + 3 \Nr{\nabla f(\ones\bfx^\star)}^2,
 \end{align}
 where ${\Q(s):=\EE{ \|\mx(s)-\bfx^\s\|^2}}$, the inequality $(a)$ follows from $(a+b+c)^2\leq 3(a^2+b^2+c^2)$, and step~$\rm{(b)}$ holds due to Assumption~\ref{asm:f}-(a). Plugging~\eqref{eq:sum-Aq} and~\eqref{eq:normgrad} into~\eqref{eq:term4-2} and~\eqref{eq:Egradient}, we arrive at
  \begin{align}\label{eq:Mtmain:exp-2}
       &\EE{ \Nr{\sum_{s=1}^{t-1} \alpha(s)  P(t:s)\nabla f(X(s))}^2}\nonumber\\
       &\leq \frac{6}{\lambda} \kap^{\frac{1}{2}} \li^2 \sum_{s=1}^{t-1}  \alpha^2(s) \pi(t:s) \lp \C(s) \!+\! \Q(s) \rp + \frac{6}{\lambda} \kap^{\frac{1}{2}} \Nr{\nabla f(\ones\bfx^\star)}^2 \sum_{s=1}^{t-1} \alpha^2(s) \pi(t:s).
   \end{align}
   To bound the second term in~\eqref{eq:Mtmain:exp-2}, we can write
 \begin{align}\label{eq:Egradient-cross} 
	 \sum_{s=1}^{t-1}  \alpha^2(s) \pi(t:s)  
       & \stackrel{\rm{(a)}}{=} \sum_{s=1}^{t-1}\alpha^2(s)\lb \beta(s)\kap^{\frac{1}{2}}\prod_{k=s+1}^{t-1}\lp 1-\lambda \beta(k)\rp^{\frac{1}{2}}\rb\nonumber\\
       &\stackrel{\rm{(b)}}{\leq} \sum_{s=1}^{t-1}\alpha^2(s)\beta(s)\kap^{\frac{1}{2}}\prod_{k=s+1}^{t-1}\lp 1-\frac{\lambda}{2} \beta(k)\rp\nonumber\\
       & = \azr^2\bzr\kap^{\frac{1}{2}} \sum_{s=1}^{t-1} \frac{1}{s^{2\nu+\mu}}\prod_{k=s+1}^{t-1}\hspace{-1pt}\lp 1\hspace{-1pt}-\hspace{-1pt} \frac{\lambda\bzr}{2}\frac{1}{k^\mu}\rp\nonumber\\
       & \stackrel{\rm{(c)}}{\leq}  \azr^2\bzr\kap^{\frac{1}{2}} A(\lambda \bzr/2, 2\nu+\mu, \mu) t^{-2\nu},
\end{align}
 for every $t\geq \Tfirst$ where 
   \begin{align}\label{eq:def:Tfirst}
       \Tfirst:=\left\lceil \lp 8\nu/ \lambda \bzr\rp^{\frac{1}{1-\mu}} \right\rceil.
   \end{align}
Note that in the chain of inequalities in~\eqref{eq:Egradient-cross}, we used the definition of $\pi(t:s)$ from \eqref{eqn:pidef} in step~\rm{(a)}, the inequality in~\rm{(b)} follows from ${\sqrt{1-x}\leq 1-x/2}$ for $x\leq 1$, and we used Lemma~\ref{lm:sum_dl_gm} with ${(\sigma, \delta, \tLM)=(2\nu+\mu, \mu, \Tfirst)}$ in~\rm{(c)}. Using~\eqref{eq:Egradient-cross} in~\eqref{eq:Mtmain:exp-2}, we get
\begin{align}\label{eq:Mtmain:exp-3}
       \EE{ \Nr{\sum_{s=1}^{t-1} \alpha(s)  P(t:s)\nabla f(X(s))}^2}\leq \frac{6}{\lambda} \kap^{\frac{1}{2}} \li^2 \sum_{s=1}^{t-1}  \alpha^2(s) \pi(t:s) \lp \C(s) \!+\! \Q(s) \rp + \cvxThree t^{-2\nu},
   \end{align}
   where $\cvxThree:=\frac{6}{\lambda} \kap \Nr{\nabla f(\ones\bfx^\star)}^2\azr^2\bzr A(\lambda \bzr/2, 2\nu\!+\!\mu, \mu) $ is a constant. Finally, using the bounds obtained in~\eqref{eq:X-mx-term} and~\eqref{eq:Mtmain:exp-3}   in~\eqref{eq:Mtmain:exp},  we get 
   \begin{align}\label{eq:conserror}
     \C(t) \leq 2\cvxTwo t^{-\mu} + \frac{12}{\lambda} \kap^{\frac{1}{2}} \li^2 \sum_{s=1}^{t-1}  \alpha^2(s) \pi(t:s) \lp \C(s) \!+\! \Q(s) \rp +2\cvxThree t^{-2\nu}.
 \end{align}
\subsection{Average State Distance to the Optimal Point}
\label{sec:avg-to-opt}
  Now, we derive an upper bound for the average distance between the mean of the agents' states and the global optimal point, i.e., ${\Q(t)=\EE{ \|\mx(t)-\bfx^\s\|^2}}$, where  $\bfx^\s$ is the minimizer of the function $f(\bfx)$. Recall that ${\mx(t) = \bfr^T X(t) =\sum_{i=1}^n r_i \bfx_i(t)}$ and ${\bfr^T W(t) = \bfr^T}$. Hence, multiplying both sides of~\eqref{eq:dyn-col-mat} by $\bfr^T$, we get
    \begin{align*}
     \mx(t+1) = \mx(t)+\et\bfr^T\ER{t} - \at\et\bfr^T\nabla f(X(t)).  
    \end{align*}

    We define ${g(t): = \bfr^T \nabla f(X(t))= \sum_{i=1}^{n}r_i\nabla f_i(\bfx_i(t))}$ and  ${\mg(t):=\nabla f(\mx(t))= \sum_{i=1}^{n}r_i\nabla f_i(\mx(t))}$.
    Hence, we can write
\begin{align}\label{eq:e_1}
    {\EE{\|\mx(t+1)-\bfx^\s\|^2\middle|\F_t}}
    &= {\EE{\|\mx(t)+\et\bfr^T\ER{t} - \at\et g(t)-\bfx^\s\|^2\middle|\F_t}}\cr
    & =\! {\|\mx(t)\!-\!\bfx^\s \!-\! \at\et g(t)\|^2}\!+\!{\EE{ \|\et\bfr^T\ER{t}\|^2\middle|\F_t}}\!,\ \ \ \ \ \ 
\end{align} 
where the last equality follows from the fact that $X(t)$ is measurable with respect to $\mathcal{F}_t$ and Assumption~\ref{assum:neighbor} implying ${\EE{\et\bfr^T\ER{t}\middle|\F_t}=0}$, which leads to
\[2\la \mx(t)- \at\et g(t)-\bfx^\s,\EE{\et\bfr^T\ER{t}\middle|\F_t}\ra=0.\]
    Taking expectations of both sides of~\eqref{eq:e_1}, and using the tower rule, we get 
    \begin{align}\label{eq:e_1_EE}
    \Q(t+1) = \EE{\|\mx(t+1)-\bfx^\s\|^2}
     =\! \EE{\|\mx(t)\!-\!\bfx^\s \!-\! \at\et g(t)\|^2}\!+\!\EE{ \|\et\bfr^T\ER{t}\|^2}\!.\ \ \ \ \ \  
    \end{align}
    In order to  bound the first term in~\eqref{eq:e_1}, we use Lemma~\ref{lm:nor_squr} for vectors  $\bfu=\mx(t)-\bfx^\s- \at\et\mg(t)$ and  ${\bfv=\alpha(t)\beta(t) (\mg(t)-g(t))}$, and parameter $\theta= \frac{\m L}{\m+L} \at\et$. Hence, we can write 
   \begin{align}\label{eq:e_3}
   \|\mx(t)-\bfx^\s-\at\et g(t)\|^2 
   &=\|\mx(t)\!-\!\bfx^\s\!-\! \at\et\mg(t)+\at\et\mg(t)\!-\!\at\et g(t)\|^2 \cr
      & \leq \lp 1+ \frac{\m L}{\m + L} \at\et\rp \|\mx(t)-\bfx^\s- \at\et\mg(t)\|^2 \cr
      &\qquad + \at\et\lp \at\et+\frac{\m + L}{\m L}\rp  \|\mg(t)-g(t)\|^2.
   \end{align}
   Next, we use Lemma~\ref{lm:sc_bnd} to bound the first term in~\eqref{eq:e_3}. Note that  Assumptions~\ref{asm:f}-(a) and~\ref{asm:f}-(b) (and Remark~\ref{rem:asm:f}) guarantee that the conditions of Lemma~\ref{lm:sc_bnd} are satisfied. Thus, we have 
   \begin{align*}
     \la \mx(t)-\bfx^\s, \nabla f(\mx(t))  \ra \geq c_1\|\nabla f(\mx(t))\|^2 +c_2 \|\mx(t)-\bfx^\s\|^2,
   \end{align*}
   or equivalently, 
   \begin{align}\label{eq:from-lm:sc_bnd}
     \la \mx(t)-\bfx^\s, \mg(t)  \ra \geq c_1\|\mg(t)\|^2 +c_2 \|\mx(t)-\bfx^\s\|^2,
   \end{align}
   where $c_1 = \frac{1}{\m + L}$ and $c_2=\frac{\m L}{\m + L}$. Therefore, for the first term in~\eqref{eq:e_3}, we can write
   \begin{align}\label{eq:e_4}
       \|\mx(t) - \bfx^\s - \at\et \mg(t)\|^2 &=\|\mx(t) - \bfx^\s\|^2  +  \alpha^2(t)\beta^2(t)\|\mg(t)\|^2 -  2\at\et \la \mx(t)-\bfx^\s, \mg(t)\ra\nonumber\\
       & \leq (1-2c_2\at\et)\|\mx(t)-\bfx^\s\|^2   +\at\et (\at\et-2c_1)\|\mg(t)\|^2.\ \ \ 
   \end{align}
   Let us set
   \begin{align}\label{eq:tau_0}
    \TStrt : = 
    \left\lceil \lp \frac{\azr\bzr}{2c_1}\rp^{\frac{1}{\mu+\nu}} \right\rceil
    =
    \left\lceil \lp \frac{\azr\bzr(\m +L)}{2}\rp^{\frac{1}{\mu+\nu}} \right\rceil,
   \end{align}
   such that $\at\et \leq 2c_1$ for any $t\geq \TStrt$. Hence, for $t\geq \TStrt$, the second term in~\eqref{eq:e_4} is non-positive, and thus 
   \begin{align*}
       \left\|\mx(t)\hspace{-1pt}-\hspace{-1pt}\bfx^\s\hspace{-3pt}-\hspace{-1pt}\at\et \mg(t)\right\|^2 \hspace{-1pt}\leq \hspace{-1pt} (1\hspace{-1pt}-\hspace{-1pt}2c_2\at\et)\|\mx(t)\hspace{-1pt}-\hspace{-1pt}\bfx^\s \|^2.  
   \end{align*}
   Taking expectations of both sides, we get
   \begin{align}\label{eq:e_5_EE}
       \EE{\left\|\mx(t)\hspace{-1pt}-\hspace{-1pt}\bfx^\s\hspace{-1pt}-\hspace{-1pt}\at\et \mg(t)\right\|^2} \leq \hspace{-1pt} (1\hspace{-1pt}-\hspace{-1pt}2c_2\at\et)\Q(t).\
   \end{align}
   The average of the second term in~\eqref{eq:e_3} can be bounded as
   \begin{align}\label{eq:g_b_g}
       \mathbb{E} \big[\|\mg(t)-g(t)\|^2\big]&= \mathbb{E}\hspace{-1pt}\left[\left\|\sum_{i=1}^{n}\hspace{-1pt}r_i \lp(\nabla f_i(\mx(t))\!-\!\nabla f_i(\bfx_i(t))\rp \right\|^2\hspace{-1pt}\right]\cr
       &\stackrel{\rm{(a)}}{\leq} \EE{\sum_{i=1}^{n}r_i\|\nabla f_i(\mx(t))-\nabla f_i(\bfx_i(t))\|^2}\nonumber\\
       &\stackrel{\rm{(b)}}{\leq} \li^2\sum_{i=1}^{n}r_i \EE{ \|\mx(t)-\bfx_i(t)\|^2}\nonumber\\
       & = \li^2 \EE{\Nr{\ones\bar{\bfx}(t)-X(t)}^2}
       = \li^2 \C(t),
       \end{align} 
       where $\rm{(a)}$ follows from the convexity of $\|\cdot\|^2$ and the inequality in~$\rm{(b)}$ holds due to Assumption~\ref{asm:f}-(a).
       
   Taking expectations of both sides of~\eqref{eq:e_3} and recalling that $c_2=\m \li/(\m+\li)$, we arrive at
  \begin{align}\label{eq:e_6}
   \mathbb{E}&\left[\|\mx(t)-\bfx^\s-\at\et g(t)\|^2\right] \nonumber\\
   & \leq \lp 1+ c_2 \at\et\rp \EE{\|\mx(t)-\bfx^\s- \at\et\mg(t)\|^2}  + \at\et\lp \at\et+1/c_2 \rp  \EE{\|\mg(t)-g(t)\|^2}\nonumber\\
       &\stackrel{\rm{(c)}}{\leq} (1+c_2\at\et )(1-2c_2\at\et)\Q(t) + \at\et (\at\et + 1/c_2) \li^2\C(t) \cr
       &\stackrel{\rm{(d)}}{\leq} (1-c_2\at\et)\Q(t) + \at\et (\at\et + 1/c_2)\li^2\C(t),
   \end{align}
   where the inequality in~${\rm (c)}$ follows from~\eqref{eq:e_5_EE} and~\eqref{eq:g_b_g},  and~${\rm (d)}$ holds since 
   \[(1+c_2\at\et )(1-2c_2\at\et)\leq 1-c_2\at\et.\]
   From Assumption~\ref{assum:neighbor}, we get
   \begin{align}\label{eq:EE}
       \lb\EE{|\ER{t}E^T(t)| \md\F_t}\rb_{ij}&=\EE{|\nt_i(t) \nt_j^T(t) |\md \F_t} \cr 
        &\leq \sqrt{\EE{\|\nt_i(t) \|^2 \md \F_t} \EE{\|\nt_j(t)\|^2 \md \F_t} }\leq \sm,
   \end{align}
    for all $1\leq i,j\leq n$. Thus, for the second term in~\eqref{eq:e_1}, we arrive at
   \begin{align}\label{eq:quant-noise}
        \EE{\left\|\bfr^T \ER{t}
        \right\|^2 \md\F_t}
         &=\bfr^T\EE{\ER{t}E^T(t)
         \md\F_t}  \bfr
        \nonumber\\
        &\leq\bfr^T\EE{|\ER{t}E^T(t)|
         \md\F_t}  \bfr \nonumber\\
         &\leq \bfr^T (\sm \ones \ones^T) \bfr = \sm.
   \end{align}
   Note that  we used the fact that $\bfr^T\ones =1$.
   Taking expectations from both sides of~\eqref{eq:quant-noise}, and using the tower rule, we arrive at
   \begin{align}\label{eq:quant-noise_EE}
       \EE{\left\|\bfr^T \ER{t}
        \right\|^2} = \EE{\EE{\left\|\bfr^T \ER{t}
        \right\|^2 \middle|\F_t}}
         \leq  \sm.
   \end{align}
   
   Using~\eqref{eq:e_6} and~\eqref{eq:quant-noise_EE} 
    in~\eqref{eq:e_1_EE}, we can write 
   \begin{align}\label{eq:e_8}
       \Q(t+1) &\leq (1-c_2\at\et)\Q(t)\!+\! \at\et (\at\et + 1/c_2)\li^2 \C(t) + \sm\beta^2(t). 
   \end{align}
   
    Now, we have two inequalities between $\C(t)$ and $\Q(t)$, namely~\eqref{eq:conserror} and~\eqref{eq:e_8}. In the following, we use these inequalities to show that both $\C(t)$ and $\Q(t)$ vanish as $t$ grows, and conclude that the state of all agents converges to $\bfx^\star$ in expectation. More precisely, we show that 
   \begin{align}\label{eq:ind-main}
        \Q(t) \leq D t^{-\min(\mu-\nu,2\nu)}, \quad \textrm{and} \quad 
        \C(t) \leq \frac{c^2_2}{2\li^2} D t^{-\min(\mu, 2\nu)},
   \end{align}
   for every $t\geq \Ts $, where
    \begin{align}\label{eq:D}
        D:=  \max \lp \max_{1\leq t\leq T_0} \Q(t) t^{\min(\mu-\nu,2\nu)}, 
        \max_{1\leq t\leq T_0} \frac{2\li^2}{c_2^2} \C(t) t^{\min(\mu,2\nu)}, \frac{8\sm \bzr}{c_2\azr}, \frac{8(\cvxTwo+\cvxThree)\li^2}{c_2^2} 
        \rp,
    \end{align}
and $\Ts=\max(\Tzero, \Tfirst, \TStrt, \Ta, \Tb, \Tc, \Td, \Te )$. Note that $\Tzero$,
$\Tfirst$, and $\TStrt$ are defined above in~\eqref{eq:def:T0},  ~\eqref{eq:def:Tfirst}, ~\eqref{eq:tau_0},
and $\Ta$, $\Tb$, $\Tc$, $\Td$ and $\Te$  will be determined as in ~\eqref{eq:Ta},~\eqref{eq:Tb},~\eqref{eq:Tc},~\eqref{eq:Td}, and~\eqref{eq:Te},    respectively.

We use induction to prove~\eqref{eq:ind-main}. First, for the starting point $\Ts$ the definition $D$ in~\eqref{eq:D} implies $\Q(\Ts) \leq D t^{-\min(\mu-\nu, 2\nu)}$ and $\C(\Ts) \leq \frac{c^2_2}{2\li^2} D t^{-\min(\mu,2\nu)}$. Next, we assume that~\eqref{eq:ind-main} holds for every $s\leq t$, and prove that $t+1$, i.e, we show that  $\Q(t\!+\!1)\leq D t^{-\min(\mu-\nu, 2\nu)}$ and $\C(t\!+\!1)\leq \frac{c^2_2}{2\li^2} D t^{-\min(\mu, 2\nu)}$.

Using the facts that $\alpha(t)\leq \alpha(\Ts)\leq \alpha(\TStrt)$ and $\beta(t)\leq \beta(\Ts)\leq \beta(\TStrt)$ for any $t\geq \Ts \geq \TStrt$, we can write
\begin{align}\label{eq:c2ab}
    c_2\at\et & \leq c_2\alpha(\TStrt)\beta(\TStrt)
     = c_2 \azr \bzr \frac{1}{T_3^{\mu+\nu}} \stackrel{\rm{(a)}}{\leq}  c_2 \azr \bzr \frac{2c_1}{\azr \bzr} =2c_2 c_1 = \frac{2\m\li}{(\m+\li)^2} \stackrel{\rm{(b)}}{\leq} \frac{1}{2},
\end{align}
where in~$\rm{(a)}$ we used the definition of $\TStrt$ in~\eqref{eq:tau_0},  and~$\rm{(b)}$ follows from the fact that $(a+b)^2\geq 4ab$. Then, from~\eqref{eq:e_8} we have
\begin{align}\label{eq:Q_t-ind}
    \Q(t+1) 
    &\leq (1-c_2\at\et)\Q(t)\!+\! \at\et (c_2 \at\et + 1) \frac{\li^2}{c_2} \C(t) + \sm\beta^2(t)\nonumber\\ 
    & \stackrel{\rm{(c)}}{\leq}  
    (1-c_2\at\et)\Q(t)\!+\! \frac{3}{2} \at\et  \frac{\li^2}{c_2} \C(t) + \sm\beta^2(t)\nonumber\\
    &\stackrel{\rm{(d)}}{\leq}
    (1-c_2\at\et)Dt^{-\min(\mu-\nu, 2\nu)}\!+\!  \frac{3}{2} \at\et  \frac{\li^2}{c_2} \frac{c_2^2}{2\li^2}D t^{-\min(\mu, 2\nu)} + \sm\beta^2(t)\nonumber\\ 
    &= \lp 1-\frac{1}{4} c_2\at\et\rp Dt^{-\min(\mu-\nu, 2\nu)}+   \sm \beta^2(t)\nonumber\\
    & = \lp 1-\frac{1}{8} c_2\at\et\rp Dt^{-\min(\mu-\nu, 2\nu)}+  \et \lp \sm \et - \frac{1}{8} c_2 D \alpha(t)  t^{-\min(\mu-\nu, 2\nu)}\rp,
\end{align}
where~$\rm{(c)}$ follows from~\eqref{eq:c2ab}, for the inequality in~$\rm{(d)}$ we used the induction assumption. Note that for the last term in~\eqref{eq:Q_t-ind} we have 
\begin{align}
\label{eq:induc-Q:1}
\frac{1}{8} c_2 D \alpha(t)  t^{-\min(\mu-\nu, 2\nu)} \!=\! \frac{c_2 D \azr}{8} t^{-(\nu+ \min(\mu-\nu, 2\nu))} =    \frac{c_2 D \azr}{8} t^{-\min(\mu, 3\nu)} \!\geq\! \frac{c_2 D \azr}{8} t^{-\mu}  \!\geq\! \sm \bzr t^{-\mu} =  \sm \beta(t),
\end{align}
where the last inequality holds provided that $D\geq \frac{8\sm \bzr}{c_2 \azr}$. Plugging~\eqref{eq:induc-Q:1} into~\eqref{eq:Q_t-ind}, we arrive at 
\begin{align}\label{eq:Q_t-ind:2}
    \Q(t+1) 
    &\leq  \lp 1-\frac{1}{8} c_2\at\et\rp D\cdot t^{-\min(\mu-\nu, 2\nu)} \stackrel{\rm{(e)}}{\leq} D \cdot (t+1)^{-\min(\mu-\nu, 2\nu)},
\end{align}
where the inequality in~{\rm (e)} follows from Lemma~\ref{lm:t->t+1} for $(a,b,c) = (c_2\azr \bzr/8, \mu+\nu, \min(\mu-\nu, 2\nu))$, and $\mu+\nu\neq 1$, which holds for 
\begin{align}\label{eq:Ta}
    t\geq \Ta := \left\{
    \begin{array}{ll}
    \left\lceil \lp \frac{8\min(\mu-\nu, 2\nu) }{c_2\azr \bzr}\rp^{\frac{1}{1-\mu-\nu}}\right\rceil & \textrm{if $\mu+\nu<1$},\\
    1 & \textrm{if $\mu+\nu=1$}.
    \end{array}
    \right.
\end{align}
Note that  when $\mu+\nu=1$, Lemma~\ref{lm:t->t+1} implies that the inequality in~\eqref{eq:Q_t-ind:2} holds for every $t\geq 1$, provided that $c_2 \azr \bzr \geq 8 \min(\mu-\nu,2\nu)$.
This completes the proof of induction for $\Q(t)$.

Next, for $\C(t)$ from~\eqref{eq:conserror}, we can write
\begin{align}\label{eq:induc:C:1}
    \C(t\!+\!1) &\leq 2\cvxTwo (t\!+\!1)^{-\mu} + \frac{12}{\lambda} \kap^{\frac{1}{2}} \li^2 \sum_{s=1}^{t}  \alpha^2(s) \pi(t\!+\!1:s) \lp \C(s) \!+\! \Q(s) \rp + 2\cvxThree (t\!+\!1)^{-2\nu}.
    \end{align}
Recall the induction assumption  $\C(s) \leq \frac{c_2^2}{2\li^2} D s^{-\min(\mu, 2\nu)}$ and $\Q(s) \leq  D s^{-\min(\mu-\nu, 2\nu)}$ for every $s\leq t$. Therefore, 
we can bound the summations in~\eqref{eq:induc:C:1} using a chain of inequalities similar to those used in~\eqref{eq:Egradient-cross} and Lemma~\ref{lm:sum_dl_gm} for $(a,\sigma, \delta)=(\lambda\bzr/2, 2\nu+\mu+\min(\mu,2\nu),\mu)$.  We have
\begin{align}\label{eq:bound-sum-C}
    \sum_{s=1}^{t}    \alpha^2(s) \pi(t\!+\!1:s) \C(s)
    &\leq  
    \sum_{s=1}^{t}\alpha^2(s)\C(s) \lb \beta(s)\kap^{\frac{1}{2}}\prod_{k=s+1}^{t}\lp 1-\lambda \beta(k)\rp^{\frac{1}{2}}\rb\nonumber\\
    &\leq \sum_{s=1}^{t}\alpha^2(s)\beta(s)\C(s) \kap^{\frac{1}{2}}\prod_{k=s+1}^{t}\lp 1-\frac{\lambda}{2} \beta(k)\rp\nonumber\\
    & =  \frac{c_2^2}{2\li^2}D \azr^2\bzr\kap^{\frac{1}{2}} \sum_{s=1}^{t} \frac{1}{s^{2\nu+\mu + \min(\mu, 2\nu)}}\prod_{k=s+1}^{t}\hspace{-1pt}\lp 1\hspace{-1pt}-\hspace{-1pt} \frac{\lambda\bzr}{2}\frac{1}{k^\mu}\rp\nonumber\\
    & \leq \frac{c_2^2}{2\li^2}D \azr^2\bzr\kap^{\frac{1}{2}} A(\lambda \bzr/2, 2\nu+\mu+\min(\mu, 2\nu), \mu) (t+1)^{-2\nu-\min(\mu, 2\nu)},
\end{align}
for every 
\begin{align}\label{eq:Tb}
t\geq \Tb:=\left\lceil \lp \frac{8\nu+4\min(\mu, 2\nu)}{\lambda \bzr}\rp ^{\frac{1}{1-\mu}} \right\rceil.
\end{align}
Similarly, we have 
 \begin{align}\label{eq:bound-sum-Q}
    \sum_{s=1}^{t}    \alpha^2(s) \pi(t\!+\!1:s) \Q(s)
    &\leq  
 \sum_{s=1}^{t}\alpha^2(s)\Q(s) \lb \beta(s)\kap^{\frac{1}{2}}\prod_{k=s+1}^{t}\lp 1-\lambda \beta(k)\rp^{\frac{1}{2}}\rb\nonumber\\
       &\leq \sum_{s=1}^{t}\alpha^2(s)\beta(s)\Q(s) \kap^{\frac{1}{2}}\prod_{k=s+1}^{t}\lp 1-\frac{\lambda}{2} \beta(k)\rp\nonumber\\
       & =  D \azr^2\bzr\kap^{\frac{1}{2}} \sum_{s=1}^{t} \frac{1}{s^{2\nu+\mu + \min(\mu-\nu, 2\nu)}}\prod_{k=s+1}^{t}\hspace{-1pt}\lp 1\hspace{-1pt}-\hspace{-1pt} \frac{\lambda\bzr}{2}\frac{1}{k^\mu}\rp\nonumber\\
       & \leq  D\azr^2\bzr\kap^{\frac{1}{2}} A(\lambda \bzr/2, 2\nu+\mu+\min(\mu-\nu, 2\nu), \mu) (t+1)^{-2\nu-\min(\mu-\nu, 2\nu)},
\end{align}
for every 
\begin{align}\label{eq:Tc}
t\geq \Tc:=\left\lceil \lp \frac{8\nu+4\min(\mu-\nu, 2\nu)}{\lambda \bzr}\rp ^{\frac{1}{1-\mu}} \right\rceil.
\end{align}
Plugging~\eqref{eq:bound-sum-C} and~\eqref{eq:bound-sum-Q} in~\eqref{eq:induc:C:1}, we arrive at
\begin{align}\label{eq:induc:C:2}
    \C(t\!+\!1) 
    & \leq 2\cvxTwo(t\!+\!1)^{-\mu} + \cvxFive D (t+1)^{-2\nu -\min(\mu,2\nu)}
    + \cvxSix D (t+1)^{-2\nu -\min(\mu-\nu,2\nu)}
     + 2\cvxThree  (t\!+\!1)^{-2\nu}\nonumber\\
     & \leq \lp 2\cvxTwo +
     \cvxFive D(t+1)^{-2\nu} + \cvxSix D (t+1)^{-\nu} + 
     2\cvxThree\rp (t+1)^{-\min(\mu, 2\nu)},
    \end{align}
    where 
    $\cvxFive:=\frac{6}{\lambda} \kap c_2^2 \azr^2\bzr A(\lambda \bzr/2, 2\nu+\mu+\min(\mu, 2\nu), \mu) $ and $\cvxSix:=\frac{12}{\lambda} \kap  \li^2 \azr^2\bzr A(\lambda \bzr/2, 2\nu+\mu+\min(\mu-\nu, 2\nu), \mu) $. Hence, in order to complete the induction for $\C(t)$, it suffices to show that 
    \begin{align}\label{eq:induc:C:3}
        2\cvxTwo +
     \cvxFive D (t+1)^{-2\nu} + \cvxSix D (t+1)^{-\nu} + 
     2\cvxThree \leq \frac{c_2^2}{2\li^2} D,
    \end{align}
    for $t\geq \Ts$. Note for 
    \begin{align}\label{eq:Td}
    t\geq \Ts \geq \Td:=\left \lceil \lp\frac{8\li^2 \cvxFive}{c_2^2 }\rp^{\frac{1}{2\nu}}\right\rceil, 
    \end{align}
    we have $\cvxFive(t+1)^{-2\nu} \leq \frac{c_2^2}{8\li^2}$. Similarly, when 
    \begin{align}\label{eq:Te}
    t\geq \Ts \geq \Te:=\left \lceil \lp\frac{8\li^2 \cvxSix}{c_2^2 }\rp^{\frac{1}{\nu}}\right\rceil, 
    \end{align}
    we get $\cvxSix(t+1)^{-\nu} \leq \frac{c_2^2}{8\li^2}$. Therefore, for $t\geq \max(\Td, \Te)$,  the inequality in~\eqref{eq:induc:C:3} holds provided that 
    \begin{align}\label{eq:induc:C:4}
        2\cvxTwo  +
     2\cvxThree \leq \frac{c_2^2}{4\li^2} D,
    \end{align}
    which clearly holds for $D\geq \frac{8(\cvxTwo + \cvxThree)\li^2}{c_2^2}$. This concludes the induction proof for $\C(t)$.

    \subsection{Total State Deviation from the Optimum Solution}
    \label{sec:prf-combine}
    In the previous section, we identified a rate (and conditions) for which the deviation between the states and their average $\C(t) =\EE{ \Nr{X(t)-\ones\bar{\bfx}(t)}^2}$ vanishes as $t$ grows. We established a similar result for the distance between the states' average and the optimum solution $\Q(t) = \EE{ \Nr{\bar{\bfx}(t)-\bfx^\s}^2}$.  Combining these bounds using Lemma~\ref{lm:nor_squr} with $\theta=1$, we can conclude the proof of Theorem~\ref{thm:strng_cvx}. In particular, we have 
    \begin{align*}
        \EE{ \Nr{X(T)\!-\!\ones\bfx^\s}^2}
        \!& = \EE{ \Nr{X(T)-\ones\bar{\bfx}(T)+\ones\bar{\bfx}(T)-\ones\bfx^\s}^2}\cr
        & \leq 2\lp \EE{ \Nr{X(T)-\ones\bar{\bfx}(T)}^2} + \EE{ \Nr{\ones\bar{\bfx}(T)-\ones\bfx^\s}^2}\rp\cr
        & = 2\lp \EE{ \Nr{X(T)-\ones\bar{\bfx}(T)}^2} + \EE{ \|\bar{\bfx}(T)-\bfx^\s\|^2}\rp\cr
        & \leq \frac{c_2^2}{\li^2} D  t^{-\min(\mu,2\nu)} \! +\! 2D t^{\min(\mu-\nu,2\nu)},
    \end{align*}
    for every ${T\geq \Ts = \max(\Tzero,\Tfirst, \TStrt,\Ta,\Tb, \Tc, \Td, \Te)}$. 
    This implies the claim of Theorem~\ref{thm:strng_cvx}, where 
    \begin{align}\label{eq:C-vals}
    \cOne := c_2^2/\li^2, \qquad\qquad \cTwo:= 2D.
    \end{align}
    Note that when $\mu+\nu=1$, the bound on $\Q(t)$ only holds when $c_2 \azr \bzr \geq 8 \min(\mu-\nu, 2\nu)$.
    This completes the proof Theorem~\ref{thm:strng_cvx}.

\section{Conclusion}
\label{sec:con}
We have studied distributed optimization over time-varying networks  suffering from noisy and imperfect sharing of information. Inspired by the original averaging-based distributed optimization algorithm with the diminishing step-size, we showed that for the class of strongly convex cost functions, including a damping mechanism for the imperfect incoming information from neighboring agents, leads to convergence to true optimizer in $L_2$ sense for various choices of the damping and diminishing step-size parameters. For the proposed algorithm, we obtain a convergence rate as a function of the damping and diminishing step-size parameters. Optimizing the resulting rate over the set of feasible parameters leads to the convergence rate  $\BO(T^{-1/2})$. 

Several avenues are left for future research, including the study of a similar mechanism with push-sum algorithm, distributed optimization over random networks, and non-convex distributed optimization with imperfect sharing of information. 

\bibliographystyle{ieeetr}
\bibliography{ref}

\appendix
\section*{Appendix: Proof of The Auxiliary Lemmas}
\label{sec:prf:aux}

In this section, we provide the proofs of auxiliary lemmas. 

	{\it Proof of Lemma~\ref{lemma:transition2}}:
	Due to the separable nature of $\Nr{\cdot}$, i.e.,  $\Nr{U}^2=\sum_{j=1}^d\Nr{U^{j}}$, without loss of generality, we may assume that $d=1$. Thus, let $U=\bfu\in \R^n$. Define $\V:\R^n\to \R^+$ by 
	    \begin{align}
	        \V(\bfu) := \Nr{\bfu-\ones \bfr^T\bfu}^2=\sum_{i=1}^n r_i(u_i-\bfr^T\bfu)^2.
	    \end{align}
	    For notational simplicity, let $\bfu(s)=\bfu=\begin{bmatrix} u_1 & u_2 & \ldots & u_n \end{bmatrix}$, and  $\bfu(k+1)=A(k)\bfu(k)$. In addition with a slight abuse of notation, we  denote $\V(\bfu(k))$ by $\V(k)$ for $k=s,\ldots,t$. 
	    
	    Using Theorem 1 in  \cite{touri2011existence}, we have 
	    \begin{align}\label{eqn:decreaseV}
	       \V(t)=\V(s)-\sum_{k=s}^{t-1}\sum_{i<j}H_{ij}(k)(u_i(k)-u_j(k))^2,
	    \end{align}
	    where $H(k)=A^T(k)\diag(\bfr)A(k)$. Note that $A(k)$ is a non-negative matrix, then ${H(k)\geq \rmin A^T(k)A(k)}$, for $k=s,\ldots,t$. Also, note that since ${A(k)=(1-\beta(k))I+\beta(k)W(k)}$, then  Assumption~\ref{asm:W}-(b) implies that the minimum non-zero elements of $A(k)$ are bounded bellow by $\minW\beta(k)$. Therefore, since $\beta(k)$ is non-increasing, on the window $k=s,\ldots,s+B$, the minimum non-zero elements of $A(k)$  for $k$ in this window are lower bounded by $\minW\beta(s+B)$. Without loss of generality, assume that the entries of $\bfu$ are sorted, i.e., $u_1\leq \ldots\leq u_n$, otherwise, we can relabel the agents (rows and columns of $A(k)$s and $\bfu$ to achieve this).  Therefore, by Lemma~8 in \cite{nedic2008distributed}, for \eqref{eqn:decreaseV}, we have 
	    \begin{align}\label{eqn:decreaseV2}
	       \V&(s+B)\nonumber\\
	       &\leq\V(s)- \rmin\sum_{k=s}^{s+B-1}\sum_{i<j}[A^T(k)A(k)]_{ij}(u_i(k)-u_j(k))^2\nonumber\\ 
	       &\leq \V(s)-\frac{\minW\rmin}{2}\beta(s+B)\sum_{\ell=1}^{n-1}(u_{\ell+1}-u_{\ell})^2.
	    \end{align}
	    
	    We may comment here that although Lemma~8 in \cite{nedic2008distributed} is written for doubly stochastic matrices, and its statement is about the decrease of $\V(\bfx)$ for the special case of  $\bfr=\frac{1}{n}\ones$, but in fact, it is a result on bounding $\sum_{k=s}^{s+B-1}\sum_{i<j}[A^T(k)A(k)]_{ij}(u_i(k)-u_j(k))^2$ for a sequence of {$B$-connected} stochastic matrices $A(k)$ in terms of the minimum non-zero entries of stochastic matrices ${A(s),\ldots,A(s+B-1)}$. 
	    
	    Next, we will show that $\sum_{\ell=1}^{n-1}(u_{\ell+1}-u_{\ell})^2\geq n^{-2}\V(\bfu)$. This argument adapts a similar argument used in the proof of Theorem~18 in \cite{nedic2008distributed} to the general $\V(\cdot)$. 
	    
	    For a $\bfv\in \R^n$ with $\V(\bfv)>0$, define the quotient 
	    \begin{align}
	        h(\bfv)=\frac{\sum_{\ell=1}^{n-1}(v_{\ell+1}-v_{\ell})^2}{\sum_{i=1}^nr_i(v_i-\bfr^T\bfv)^2}=\frac{\sum_{\ell=1}^{n-1}(v_{\ell+1}-v_{\ell})^2}{\V(\bfv)}. 
	    \end{align}
	    Note that $h(\bfv)$ is invariant under scaling and translations by all-one vector, i.e., $h(\omega \bfv )=h(\bfv)$ for all non-zero $\omega\in \R$ and $h(\bfv+\omega \ones)=h(\bfv)$ for all $\omega\in \R$. Therefore, 
	    \begin{align}\label{eqn:hv}
	        \min_{\substack{v_1\leq v_2\leq \cdots\leq v_n\\\V(\bfv)\not=0}}h(\bfv)&=\min_{\substack{v_1\leq v_2\leq \cdots\leq v_n\\ \bfr^T\bfv=0,\V(\bfv)=1}}h(\bfv)\cr 
	        &=\min_{\substack{v_1\leq  v_2\leq \cdots\leq v_n\\\bfr^T\bfv=0,\V(\bfv)=1}}\sum_{\ell=1}^{n-1}(v_{\ell+1}-v_{\ell})^2. 
	    \end{align}
	    
	    Since $\bfr$ is a stochastic vector, then for a vector $\bfv$ with  ${v_1\leq \ldots\leq v_n}$ and $\bfr^T\bfv=0$, we would have ${v_1\leq \bfr^T\bfv=0\leq v_n}$. On the other hand, the fact that ${\V(\bfv)=\sum_{i=1}^nr_iv^2_i=1}$ would imply ${\max(|v_1|,|v_n|)\geq \frac{1}{\sqrt{n}}}$. Let us consider the difference sequence  $\hat{v}_\ell=v_{\ell+1}-v_{\ell}$  for $\ell=1,\ldots, n-1$, for which we have $\sum_{i=1}^{n-1} \hat{v}_\ell = v_n - v_1 \geq v_n \geq \frac{1}{\sqrt{n}}$. Therefore, the optimization problem \eqref{eqn:hv} can be  rewritten as
        \begin{align}\label{eqn:hv2}
	        \min_{\substack{v_1\leq v_2\leq \cdots\leq v_n\\\V(\bfv)\not=0}}h(\bfv)&=\min_{\substack{v_1\leq v_2\leq \cdots\leq v_n\\\bfr^T\bfv=0,\V(\bfv)=1}}\sum_{\ell=1}^{n-1}(v_{\ell+1}-v_{\ell})^2\cr 
	        &\geq 
	        \min_{\substack{\hat{v}_1,\ldots,\hat{v}_{n-1}\geq 0\\\sum_{i=1}^{n-1}\hat{v}_i\geq \frac{1}{\sqrt{n}}}}\sum_{\ell=1}^{n-1}\hat{v}_\ell^2. 
	    \end{align}
	    Using the Cauchy-Schwarz inequality, we get 
	    ${\lp \sum_{\ell=1}^{n-1} \hat{v}_\ell^2 \rp \cdot \lp \sum_{\ell=1}^{n-1} 1^2\rp \geq \big( \sum_{\ell=1}^{n-1} \hat{v}_\ell \big)^2 \geq \big( \frac{1}{\sqrt{n}} \big)^2 = \frac{1}{n}}$. Hence, 
	    
	    \begin{align}
	         \min_{\substack{v_1\leq v_2\leq \ldots\leq v_n\\
	         \V(\bfv)\not=0}}h(\bfv)\geq \frac{1}{n(n-1)}\geq \frac{1}{n^2}.
	    \end{align}
        Thus, for $v_1\leq  \ldots\leq v_n$, we have ${\sum_{\ell=1}^{n-1}(v_{\ell+1}-v_{\ell})^2\geq n^{-2}\V(\bfv)}$ (note that this inequality also holds for $\bfv\in\R^n$ with $\V(\bfv)=0$). 
	    Using this fact in \eqref{eqn:decreaseV2} implies  
	    \begin{align}\label{eqn:Vdot}
	       \V(s+B) \leq \lp 1-\frac{\minW\rmin}{2n^2}\beta(s+B)\rp\V(s).
	    \end{align}
	    
	    Applying~\eqref{eqn:Vdot} for $\Del:=\lfloor \frac{t-1-s}{B}\rfloor$ steps recursively, we get
	    \begin{align*}
	        \V(s+\Del B) \leq \prod_{j=1}^{\Del} \lp 1-\frac{\minW\rmin}{2n^2}\beta(s+jB)\rp\V(s)
	    \end{align*}
	    
Using the fact that  $(1-x)^{1/B} \leq 1-x/B$ and since $\{\beta(k)\}$ is a non-increasing sequence, we have
\begin{align*}
	1-\frac{\minW\rmin}{2n^2}\beta(s+jB) 
	&= \prod_{\ell=1}^{B} \lp 1-\frac{\minW\rmin}{2n^2}\beta(s+jB)\rp^{1/B}\\  
	&\leq 
	\prod_{\ell=1}^{B}  \lp 1-\frac{\minW\rmin}{2Bn^2}\beta(s+jB)\rp\nonumber\\
	&\leq 
	\prod_{\ell=1}^{B}  \lp 1-\lambda \beta(s+jB+\ell)\rp.
\end{align*}
Recall from~\eqref{eqn:decreaseV} that $\V(t)$ is a non-increasing function of $t$. Therefore, for ${s+\Del B \leq t -1 < s+(\Del+1) B }$ we have
\begin{align}\label{eq:B-jumps}
	  \V(t-1) &\leq \V(s+\Del B)\nonumber\\
	  &\leq \prod_{j=1}^{\Del} \lp 
	  1-\frac{\minW\rmin}{2n^2}\beta(s+jB) 
	  \rp\V(s) \nonumber\\
	  &\leq \prod_{j=1}^{\Del} 	\prod_{\ell=1}^{B}  \lp 1- \lambda \beta(s+jB+\ell)\rp \V(s) \nonumber\\
	  &= \prod_{k=s+B+1}^{s+(\Del+1)B} \lp 1- \lambda \beta(k)\rp \V(s) \nonumber\\
	  &\leq \prod_{k=s+B+1}^{t-1} \lp 1-\lambda \beta(k)\rp \V(s). 
\end{align}
	   Next, noting that $\{\beta(k)\}$ is a non-increasing sequence, we have ${\beta(k) \leq \beta(1) = \bzr}$. Thus, 
	   \begin{align}\label{eq:initial-steps}
	       \prod_{k=s+1}^{s+B} \lp 1- \lambda \beta(k)\rp &\geq \prod_{k=s+1}^{s+B} \lp 1- \lambda \bzr \rp \nonumber\\
	       &= \lp 1- \lambda \bzr \rp^B \geq 1- B\lambda \bzr. 
	   \end{align}
	   Therefore, combining~\eqref{eq:B-jumps} and~\eqref{eq:initial-steps}, we get 
	   \begin{align}\label{eq:s-to-t}
	       \V(t-1) &\leq \prod_{k=s+B+1}^{t-1} \lp 1- \lambda \beta(k)\rp  \V(s) \nonumber\\
	       &\leq 
	       \frac{\prod_{k=s+1}^{s+B} \lp 1- \lambda \beta(k)\rp }{1-B\lambda \bzr}\!\! \prod_{k=s+B+1}^{t-1} \!\!\lp 1- \lambda \beta(k)\rp  \V(s) \nonumber\\
	       &= 
	       \frac{1}{1-B\lambda \bzr} \prod_{k=s+1}^{t-1} (1-\lambda \beta(k)) \V(s). 
	   \end{align}
	   Now, we define 
	   \[
	   \Phi(t\hspace{-1pt}:\hspace{-1pt}s)=A(t\hspace{-1pt}-\hspace{-1pt}1)\cdots A(s\hspace{-1pt}+\hspace{-1pt}1)\] for ${t\geq s}$ with ${\Phi(t\hspace{-1pt}:\hspace{-1pt}t\hspace{-1pt}-\hspace{-1pt}1)=I}$. Note that   {Assumption~\ref{asm:W}-(a)} and the fact that ${A(k) = (1-\beta(k))I + \beta(k) W(k)}$ imply ${\bfr^T \Phi(t:s)  =\bfr^T}$. Then, setting $\bfu(s) = \bfu = U$ and ${\bfu(t-1) = \Phi(t:s) \bfu(s) =\Phi(t:s) U}$, we can write
	   \begin{align*}
	       \lp \Phi(t:s) -\ones \bfr^T\rp U & = \Phi(t:s)U- \ones \bfr^T U \nonumber\\
	       &=\Phi(t:s)U-\ones \bfr^T \Phi(t:s) U \nonumber\\
	       &= \bfu(t-1) - \ones \bfr^T \bfu(t-1).
	   \end{align*}
	   Therefore, using~\eqref{eq:s-to-t} we have
	   \begin{align*}
	       \big\| \big( \Phi(t:s) &-\ones \bfr^T\big) U \big\|_\bfr^2  \nonumber\\
	       &= \Nr{\bfu(t-1) - \ones \bfr^T \bfu(t-1)}^2 = \V(t-1) \cr 
	       & \leq \frac{1}{1-B\lambda \bzr} \prod_{k=s+1}^{t-1} (1-\lambda \beta(k)) \V(s)  \\
	       &= \frac{1}{1-B\lambda \bzr} \prod_{k=s+1}^{t-1} (1-\lambda \beta(k)) \Nr{\bfu - \ones \bfr^T \bfu}^2 
	       \cr
	       & \leq \frac{1}{1-B\lambda \bzr} \prod_{k=s+1}^{t-1} (1-\lambda \beta(k)) \Nr{\bfu}^2 \\
	       &=\frac{1}{1-B\lambda \bzr} \prod_{k=s+1}^{t-1} (1-\lambda \beta(k)) \Nr{U}^2, 
	   \end{align*}
	   where the second inequality follows from the fact that $\Nr{\bfu - \ones \bfr^T \bfu}^2 + \Nr{\ones \bfr^T \bfu}^2 = \Nr{\bfu}^2$. This completes the proof of the lemma. 
	\hfill $\blacksquare$

{\it Proof of Lemma~\ref{lemma:matrixprop}:} 
	  Recall that we denote the $i$th row of $A$ by $A_i$, and the $j$th column of $B$ by $B^j$. Then, applying the Cauchy-Schwartz inequality to vectors $A_i$ and $B^j$, we have ${|[AB]_{ij}|=|\la A_i,B^j\ra| \leq \lnr A_i\rnr \lnr B^j\rnr}$. Therefore, 
	  \begin{align*}
	   \lnr[AB]_i\rnr^2
	   & =\sum_{j=1}^m |[AB]_{ij}| = 
	   \sum_{j=1}^m|\la A_i,B^j\ra|^2\cr
	   & \leq \lnr A_i\rnr^2 \sum_{j=1}^m\lnr B^j\rnr^2\leq \lnr A_i\rnr^2 \lnr B\rnr_F^2.
	  \end{align*}
	  
	 Using this inequality and the definition of $\bfr$-norm, we get
	  \begin{align*}
	     \lnr AB \rnr_\bfr^2 &=  \sum_{i=1}^nr_i\lnr[AB]_i\rnr^2\leq \sum_{i=1}^nr_i\lnr A_i\rnr^2 \lnr B\rnr_F^2 = \lnr A\rnr_\bfr^2 \lnr B\rnr_F^2,
 	  \end{align*}
 	  as claimed in the lemma.
\hfill $\blacksquare$ 	  
 	  
{\it Proof of Lemma~\ref{lm:nor_squr}:} 
For two vectors $\bfu$ and $\bfv$ and any scalar $\theta>0$, we have
   \begin{align*}
       \|\bfu+\bfv\|^2 & = \|\bfu\|^2 + \|\bfv\|^2 + 2\la\bfu,\bfv \ra \cr
       &\stackrel{\rm{(a)}}{\leq}  \|\bfu\|^2 + \|\bfv\|^2 + 2 \|\bfu \| \|\bfv\| \\
       &= \|\bfu\|^2 + \|\bfv\|^2 + 2 \lp \sqrt{\theta}\|\bfu \| \cdot \frac{1}{\sqrt{\theta}}\|\bfv\|\rp \\
        &\stackrel{\rm{(b)}}{\leq}  \|\bfu\|^2 + \|\bfv\|^2 +  \theta \|\bfu \|^2 + \frac{1}{\theta}\|\bfv\|^2 \\
        &= (1+\theta)  \|\bfu\|^2 + \lp 1+ \frac{1}{\theta} \rp \|\bfv\|^2,
   \end{align*}
   where $\rm{(a)}$ follows from the Cauchy–Schwarz inequality and $\rm{(b)}$ is concluded from the geometric-arithmetic inequality. Similarly, recalling $\bfr$-norm, for matrices $U$ and $V$, we have 
   \begin{align*}
       \|U+V\|_\bfr^2 &= \sum_{i=1}^n r_i \|U_i+V_i\|^2 \cr
       &\leq \sum_{i=1}^n r_i \lb (1+\theta) \|U_i\|^2 + \lp 1+\frac{1}{\theta}\rp \|V_i\|^2\rb \cr
       &= (1+\theta) \sum_{i=1}^n r_i \|U_i\|^2 + \lp 1+\frac{1}{\theta}\rp \sum_{i=1}^n r_i \|V_i\|^2\cr
       &= (1+\theta) \|U\|^2_\bfr + \lp 1+\frac{1}{\theta}\rp  \|V\|^2_\bfr.
   \end{align*}
   This completes the proof of the lemma. 
   \hfill $\blacksquare$
   
{\it Proof of Lemma~\ref{lm:delta}:} 
    Since $\log(1-x)\leq -x$ and $\frac{a}{\tau^{\dl}}$ is a decreasing function of $\tau$, for $0<\dl<1$ we have
    \begin{align*}
        \log \prod_{k=s}^{t-1} \lp 1-\frac{a}{k^{\dl}} \rp  &= \sum_{k=s}^{t-1} \log\lp 1-\frac{a}{k^{\dl}} \rp \\
        &\leq -\sum_{k=s}^{t-1} \frac{a}{k^{\dl}}\leq - \int_{s}^{t} \frac{a}{\tau^{\dl}} d\tau\\
        &= -\frac{a}{1-\dl} \lp t^{1-\dl} - s^{1-\dl}\rp.
    \end{align*}
        Thus, we have ${\prod_{k=s}^{t-1} \lp 1-\frac{a}{k^{\dl}} \rp \leq \exp\lp -\frac{a}{1-\dl} \lp t^{1-\dl} - s^{1-\dl}\rp \rp}$. Using a similar argument for $\dl=1$, we can write
    \begin{align*}
         \log \prod_{k=s}^{t-1} \lp 1-\frac{a}{k} \rp  &= \sum_{k=s}^{t-1} \log\lp 1-\frac{a}{k} \rp \\
        &\leq -\sum_{k=s}^{t-1} \frac{a}{k}\leq - \int_{s}^{t} \frac{a}{\tau} d\tau \cr
        &= -a  \ln\lp\frac{t}{s} \rp.
     \end{align*}
    This implies that $\prod_{k=s}^{t-1} \lp 1-\frac{a}{k} \rp \leq \exp\lp -a \ln \lp \frac{t}{s}\rp \rp = \lp \frac{t}{s}\rp^{-a}$. This completes the proof of the lemma. \hfill $\blacksquare$
    
{\it Proof of Lemma~\ref{lemma:Psum}:} 
        In order to prove~\eqref{eq:esumprod}, we define ${p=\sum_{s=1}^{t-1} \beta(s)\prod_{k=s+1}^{t-1}(1-\lambda \beta(k))}$. Then, we have
        \begin{align*}
            \lambda p &= \sum_{s=1}^{t-1} \lambda \beta(s)\prod_{k=s+1}^{t-1}(1-\lambda \beta(k))\cr
                      &= \sum_{s=1}^{t-1} (1-(1-\lambda \beta(s)))\prod_{k=s+1}^{t-1}(1-\lambda \beta(k))\cr 
                      &= \sum_{s=1}^{t-1}\left[\prod_{k=s+1}^{t-1}\lp 1-\lambda \beta(k)\rp -\prod_{k=s}^{t-1}\lp1-\lambda \beta(k))\rp\right].
        \end{align*}
        Noticing that the last sum is a telescopic sum implies
        \begin{align*}
            \lambda p &= \prod_{k=t}^{t-1}\lp1-\lambda \beta(k)\rp -\prod_{k=1}^{t-1}\lp1-\lambda \beta(k)\rp\cr
            & = 1 - \prod_{k=1}^{t-1}\lp1-\lambda \beta(k)\rp.
        \end{align*}
        Dividing both sides by $\lambda\not=0$ arrives at~\eqref{eq:esumprod}.\hfill $\blacksquare$

{\it Proof of Lemma~\ref{lm:sum_dl_gm}:} 
    First, note from Lemma~\ref{lm:delta} that 
    \begin{align}
        \prod_{k=s}^{t-1}\lp 1-\frac{a}{k^\dl}\rp \leq  \exp\lp -\frac{a}{1-\dl} \lp t^{1-\dl} - s^{1-\dl}\rp \rp. 
        \end{align}
    Therefore, we get
    \begin{align}
     &\sum_{s=1}^{t-1}  \lb \frac{1}{s^{\sigma}} \prod_{k=s+1}^{t-1} \lp 1-\frac{a}{k^{\dl}} \rp \rb \cr
     &= \sum_{s=2}^{t}  \lb \lp\frac{s}{s-1}\rp^{\sigma}  \frac{1}{s^{\sigma}} \prod_{k=s}^{t-1} \lp 1-\frac{a}{k^{\dl}} \rp  \rb \nonumber\\
     &\leq 2^\sigma \lb \sum_{s=2}^{t-1}  \lb  \frac{1}{s^{\sigma}}  \prod_{k=s}^{t-1} \lp 1-\frac{a}{k^\dl}\rp\rb + \frac{1}{t^{\sigma}} \rb \nonumber\\
     &\leq 2^\sigma \lb \sum_{s=2}^{t-1}  s^{-\sigma}  \exp\lp -\frac{a}{1-\dl} \lp t^{1-\dl} - s^{1-\dl}\rp \rp + t^{-\sigma}\rb 
     \nonumber\\
     &= 2^\sigma \lb e^{ -\frac{a}{1-\dl}  t^{1-\dl}}  \sum_{s=2}^{t-1}  s^{-\sigma} \exp\lp \frac{a}{1-\dl}  s^{1-\dl} \rp + t^{-\sigma}\rb. 
     \label{eq:pr:lm:2nd}
     \end{align}
    Now, consider function  $h(\tau) = \tau^{-\sigma} \exp\lp \frac{a}{1-\dl}  \tau^{1-\dl} \rp$, with 
    \begin{align}
    \frac{dh(\tau)}{d\tau} = \lp a\tau^{-\sigma}\tau ^{-\dl} - \sigma \tau^{-\sigma-1} \rp \exp\lp \frac{a}{1-\dl}  \tau^{1-\dl} \rp.
    \end{align}
    Let $t_0 = \lceil (\sigma/a)^{\frac{1}{1-\dl}} \rceil$. Then the function $h(\tau)$ is a decreasing function for $\tau\leq t_0-1$ and  an increasing function for 
    $\tau\geq t_0$. Therefore, we can write
     \begin{align}
      \sum_{s=2}^{t_0-1}   s^{-\sigma} \exp\lp \frac{a}{1-\dl}  s^{1-\dl} \rp
    &\leq 
    \int_{1}^{t_0-1}  \tau^{-\sigma}  e^{\frac{a}{1-\dl}  \tau^{1-\dl} } d\tau, \label{eq:left-tail}\\
     \sum_{s=t_0}^{t-1}   s^{-\sigma} \exp\lp \frac{a}{1-\dl}  s^{1-\dl} \rp
    &\leq 
    \int_{t_0}^{t}  \tau^{-\sigma}  e^{\frac{a}{1-\dl}  \tau^{1-\dl}} d\tau. \label{eq:right-tail}
    \end{align}
    Summing up inequalities in~\eqref{eq:left-tail} and~\eqref{eq:right-tail}, we arrive at 
    \begin{align}
    \sum_{s=2}^{t-1}   s^{-\sigma}  \exp\lp \frac{a}{1-\dl}  s^{1-\dl} \rp
    &\leq 
    \int_{1}^{t}  \tau^{-\sigma}  e^{\frac{a}{1-\dl}  \tau^{1-\dl}} d\tau.  
    \end{align}
    Substituting this into~\eqref{eq:pr:lm:2nd}, we get 
    \begin{align}
     \sum_{s=2}^{t-1}  \lb \frac{1}{s^{\sigma}} \prod_{k=s}^{t-1} \lp 1-\frac{a}{k^{\dl}} \rp \rb 
     &\leq 
      \frac{
     \int_{1}^{t}  \tau^{-\sigma}  \exp\lp \frac{a}{1-\dl}  \tau^{1-\dl} \rp d\tau
     }
     {\exp\lp \frac{a}{1-\dl}  t^{1-\dl}\rp
     }.
     \label{eq:lm:term2:1}
     \end{align}
     Now, let us define ${p(t) := \int_{1}^{t}  \tau^{-\sigma}  \exp\lp \frac{a}{1-\dl}  \tau^{1-\dl} \rp d\tau}$ and ${q(t):=C t^{-(\sigma-\dl)}\exp\lp \frac{a}{1-\dl}  t^{1-\dl}\rp}$ for some constant $C>0$ (independent of $t$).  In the following, we will show that $p(t) \leq q(t)$ for $t \geq \TStrtTwo := \max\lc 1, \lp \frac{2(\sigma-\delta)}{a}\rp^{\frac{1}{1-\delta}} \rc$ and a proper choice of $C$. To this end, we show that $p'(\TStrtTwo)\leq q'(\TStrtTwo)$ and $p(\TStrtTwo)\leq q(\TStrtTwo)$. 
    First note that 
     \begin{align*}
         p'(t) &= t^{-\sigma}  \exp\lp \frac{a}{1-\dl}  t^{1-\dl} \rp,\\
       q'(t)  &=  C\left( a t^{-\sigma} - (\sigma - \delta)t^{-(\sigma -\dl)-1}  \right)\exp\lp \frac{a}{1-\dl}  t^{1-\dl} \rp\cr
       & =  aC\lp 1 - \frac{(\sigma-\delta)}{a} t^{-(1-\delta)} \rp t^{-\sigma}  \exp\lp \frac{a}{1-\dl}  t^{1-\dl} \rp .
     \end{align*}
     Hence, for $t \geq \TStrtTwo$ we have 
     \begin{align*}
      q'(t)  &=   aC \lp 1 - \frac{\sigma-\delta}{a} t^{-(1-\delta)} \rp t^{-\sigma}  \exp\lp \frac{a}{1-\dl}  t^{1-\dl} \rp\\
      &\geq aC \lp 1 - \frac{\sigma-\delta}{a} \left(\frac{a}{2(\sigma-\delta)}\right)  \rp    t^{-\sigma} \exp\lp \frac{a}{1-\dl}  t^{1-\dl} \rp\\
      &= \frac{a}{2}C  t^{-\sigma}  \exp\lp \frac{a}{1-\dl}  t^{1-\dl} \rp, 
     \end{align*}
     which is greater than or equal to $p'(t)$ provided that $C\geq \frac{2}{a}$. 
     It only remains to determine $C$ such that show that $p(\TStrtTwo)\leq q(\TStrtTwo)$. 
    First note that if $\TStrtTwo=1$, then $p(\TStrtTwo)=0\leq q(\TStrtTwo)$ for any $C\geq 0$. We will prove the claim for $\sigma >1$, $\sigma=1$, and $0< \sigma <1$, separately, for the case of   $\TStrtTwo :=  \lp \frac{2(\sigma-\delta)}{a}\rp^{\frac{1}{1-\delta}}$. 
     When $\sigma>1$, we  have 
     \begin{align*}
         p(\TStrtTwo) & = \int_{1}^{\TStrtTwo}  \tau^{-\sigma}  e^{ \frac{a}{1-\dl}  \tau^{1-\dl}} d\tau 
         \leq 
         \exp\lp \frac{a}{1-\dl} \TStrtTwo^{1-\dl}\rp \int_{1}^{\TStrtTwo}  \tau^{-\sigma} d\tau  \cr
         &= \frac{1-\TStrtTwo^{1-\sigma}}{\sigma-1}
         \exp\lp \frac{a}{1-\dl} \TStrtTwo^{1-\dl} \rp\cr
         &<\frac{1}{\sigma-1} 
         \exp\lp \frac{a}{1-\dl} \TStrtTwo^{1-\dl} \rp\cr
         & = \frac{1}{\sigma-1} \TStrtTwo^{\sigma-\dl} \TStrtTwo^{-(\sigma-\delta)} 
         \exp\lp \frac{a}{1-\dl} \TStrtTwo^{1-\dl} \rp\cr
         & = \frac{1}{\sigma-1}\lp\frac{2(\sigma-\dl)}{a}\rp^{\frac{\sigma-\dl}{1-\dl}} \TStrtTwo^{-(\sigma-\dl)}
         \exp\lp \frac{a}{1-\dl} \TStrtTwo^{1-\dl} \rp\cr
         &\leq q(\TStrtTwo),
      \end{align*}
      where the last inequality holds for $C\geq \frac{1}{\sigma-1}\lp\frac{2(\sigma-\dl)}{a}\rp^{\frac{\sigma-\dl}{1-\dl}}$. 
     For the case of $\sigma=1$, we have
     \begin{align*}
         p(\TStrtTwo) & = \int_{1}^{\TStrtTwo}  \tau^{-1}  \exp\lp \frac{a}{1-\dl}  \tau^{1-\dl} \rp d\tau\cr
         &\leq 
         \exp\lp \frac{a}{1-\dl} \TStrtTwo^{1-\dl} \rp \int_{1}^{\TStrtTwo}  \tau^{-1} d\tau  \cr
         &= \ln (\TStrtTwo)
         \exp\lp \frac{a}{1-\dl} \TStrtTwo^{1-\dl} \rp\cr
         & = \ln (\TStrtTwo)\TStrtTwo^{1-\dl} \TStrtTwo^{-(1-\dl)}\exp\lp \frac{a}{1-\dl} \TStrtTwo^{1-\dl} \rp \cr
         & = \frac{2}{a}\ln\lp\frac{2(1-\dl)}{a}\rp \TStrtTwo^{-(1-\dl)}\exp\lp \frac{a}{1-\dl} \TStrtTwo^{1-\dl} \rp \leq q(\TStrtTwo),
      \end{align*}
      where the last inequality holds provided that ${C\geq \frac{2}{a}\ln\lp\frac{2(1-\dl)}{a}\rp}$. 
     Lastly, for the case of $0<\sigma<1$, we can write 
     \begin{align*}
         p(\TStrtTwo) & = \int_{1}^{\TStrtTwo}  \tau^{-\sigma}  \exp\lp \frac{a}{1-\dl}  \tau^{1-\dl} \rp d\tau \cr
         & \leq 
         \exp\lp \frac{a}{1-\dl} \TStrtTwo^{1-\dl} \rp \int_{1}^{\TStrtTwo}  \tau^{-\sigma} d\tau  \cr
         &= \frac{\TStrtTwo^{1-\sigma}-1}{1-\sigma}
         \exp\lp \frac{a}{1-\dl} \TStrtTwo^{1-\dl} \rp\cr
         &<\frac{\TStrtTwo^{1-\dl}}{1-\sigma} \TStrtTwo^{-(\sigma-\dl)}
         \exp\lp \frac{a}{1-\dl} \TStrtTwo^{1-\dl} \rp\cr
         & = \frac{2(\sigma-\dl)}{a(1-\sigma)} \TStrtTwo^{-(\sigma-\dl)}
         \exp\lp \frac{a}{1-\dl} \TStrtTwo^{1-\dl} \rp \leq q(\TStrtTwo),
      \end{align*}
      where the last inequality holds for $C\geq \frac{2(\sigma-\dl)}{a(1-\sigma)}$. 
      Therefore, we have $p(t) \leq q(t)$ for $t\geq \TStrtTwo $ where $C$ is given by  
      \begin{align*}
          C= \begin{cases}
          \max\lc \frac{2}{a}, \frac{1}{\sigma-1}\lp\frac{2(\sigma-\dl)}{a}\rp^{\frac{\sigma-\dl}{1-\dl}}\rc  & \textrm{if $\sigma>1$},\\
          \max\lc \frac{2}{a}, \frac{2}{a}\ln\lp\frac{2(1-\dl)}{a}\rp \rc  & \textrm{if $\sigma=1$},\\
          \max\lc \frac{2}{a}, \frac{2(\sigma-\dl)}{a(1-\sigma)}\rc & \textrm{if $0<\sigma<1$}. 
          \end{cases}
      \end{align*}
      Plugging this into~\eqref{eq:lm:term2:1}, we get
     \begin{align*}
         \sum_{s=2}^{t-1}  \lb \frac{1}{s^{\sigma}} \prod_{k=s}^{t-1} \lp 1-\frac{a}{k^{\dl}} \rp \rb 
      &\leq    \frac{
     \int_{1}^{t}  \tau^{-\sigma}  \exp\lp \frac{a}{1-\dl}  \tau^{1-\dl} \rp d\tau
     }
     {\exp\lp \frac{a}{1-\dl}  t^{1-\dl}\rp
     }
    \cr
    & =  \frac{p(t)}{\exp\lp \frac{a}{1-\dl}  t^{1-\dl}\rp}\\
     &\leq 
     \frac{q(t)}{\exp\lp \frac{a}{1-\dl}  t^{1-\dl}\rp}
     = C t^{-(\sigma-\dl)},
     \end{align*}
     for $t\geq \TStrtTwo$. Then, continuing from~\eqref{eq:pr:lm:2nd}, we get 
     \begin{align}
     &\sum_{s=1}^{t-1}  \lb \frac{1}{s^{\sigma}} \prod_{k=s+1}^{t-1} \lp 1-\frac{a}{k^{\dl}} \rp \rb \cr
     &\leq 2^\sigma \lb \sum_{s=2}^{t-1}  \lb  \frac{1}{s^{\sigma}}  \prod_{k=s}^{t-1} \lp 1-\frac{a}{k^{\dl}} \rp  \rb + \frac{1}{t^{\sigma}} \rb \nonumber\\
     &\leq 2^\sigma \lb C t^{-(\sigma-\dl)} + t^{-\sigma} \rb\nonumber\\
     &\leq 2^\sigma \lb C + t^{-\dl} \rb t^{-(\sigma-\dl)} \nonumber\\
     &\leq 2^\sigma [ C +1] t^{-(\sigma-\dl)}.
     \end{align}
     Therefore, $A(a,\sigma, \delta)$ in the statement of the lemma is determined to be $2^{\sigma}(C+1)$. 
     \\     
     Finally, let us consider the case of  $\dl = 1$. Similar to~\eqref{eq:pr:lm:2nd}, and using  Lemma~\ref{lm:delta} for $a-\sigma +1  \neq 0$, we have 
     \begin{align*}
     \sum_{s=1}^{t-1} \hspace{-1pt} \lb \frac{1}{s^{\sigma}} \hspace{-1pt}\prod_{k=s+1}^{t-1} \lp 1-\frac{a}{k} \rp \hspace{-1pt}\rb \hspace{-1pt}
     &\leq  2^\sigma  \lb \sum_{s=2}^{t-1}  \lb \frac{1}{s^{\sigma}} \prod_{k=s}^{t-1} \lp 1-\frac{a}{k} \rp \rb +\frac{1}{t^\sigma}\rb\nonumber\\
     & \leq 2^\sigma \lb\sum_{s=2}^{t-1} \frac{1}{s^{\sigma}} \lp\frac{t}{s}\rp^{-a} +\frac{1}{t^{\sigma}}\rb \cr
     & = 2^\sigma t^{-a}\sum_{s=2}^{t-1} s^{a-\sigma} +2^\sigma t^{-\sigma}\\
     & \leq 2^\sigma  t^{-a} \int_{1}^{t} \tau^{a-\sigma} d\tau +2^\sigma t^{-\sigma}\cr
     &=2^{\sigma}t^{-a} \left| \frac{t^{a-\sigma+1} -1 }{a-\sigma+1}\right|+2^\sigma t^{-\sigma}\cr
     &= 2^{\sigma}\left| \frac{t^{-(\sigma-1)} -t^{-a} }{a-\sigma+1}\right| +2^\sigma t^{-\sigma}\cr
     &\leq 2^{\sigma} \hspace{-1pt}\lp \frac{1}{|a-\sigma+1|}+\hspace{-1pt}1\rp t^{-\min(\sigma-1,a) }.
     \end{align*} 
     This completes the proof of the lemma. \hfill $\blacksquare$ 
     
     {\it Proof of Lemma~\ref{lm:t->t+1}}:
    First let $b\neq 1$. Then, the claimed inequality is equivalent to 
    \[
    \frac{a}{t^{c+b}}\geq \frac{1}{t^c} - \frac{1}{(t+1)^{c}}.
    \]
    The mean value theorem for $h(t)=1/t^c$ implies that 
    \[
    \frac{h(t+1) - h(t)}{(t+1)-t} = h'(z) = -\frac{c}{z^{c+1}}
    \]
    for some $z\in (t,t+1)$. Therefore, we have 
    \begin{align}\label{eq:pf:lm:t->t+1}
        \frac{1}{t^c} - \frac{1}{(t+1)^c} = \frac{c}{z^{c+1}} \leq \frac{c}{t^{c+1}} = \lp \frac{c/a}{t^{1-b}}\rp \frac{a}{t^{c+b}} \leq  \frac{a}{t^{c+b}},
    \end{align}
    where the first inequality follows from $z\geq t$ and the second inequality holds for $t\geq \lp \frac{c}{a}\rp^{\frac{1}{1-b}}$. 
    
    Finally, if $b=1$ and $a\geq c$, then we have $c/t^{c+1}\leq a/t^{c+b}$, and the inequality in~\eqref{eq:pf:lm:t->t+1} holds for all values of $t\geq 1$. 
    This completes the proof of the lemma. 
    \hfill $\blacksquare$
\end{document}